\documentclass[12pt,letterpaper]{article}

\usepackage{arxiv}

\usepackage[utf8]{inputenc} % allow utf-8 input
\usepackage[T1]{fontenc}    % use 8-bit T1 fonts
\usepackage{hyperref}       % hyperlinks
\usepackage{url}            % simple URL typesetting
\usepackage{booktabs}       % professional-quality tables
\usepackage{amsfonts}       % blackboard math symbols
\usepackage{nicefrac}       % compact symbols for 1/2, etc.
\usepackage{microtype}      % microtypography
\usepackage{lipsum}		% Can be removed after putting your text content
\usepackage{graphicx}
\usepackage{natbib}
\usepackage{doi}
\usepackage{amsmath} 
\usepackage{multirow} 
\usepackage{algorithm}
\usepackage{algpseudocode} 
\usepackage{soul,xcolor} 
\usepackage{subcaption} 
\usepackage{graphicx}

\usepackage{tikz}
\usetikzlibrary{matrix} % for block alignment
\usetikzlibrary{arrows} % for arrow heads
\usetikzlibrary{arrows.meta}
\usetikzlibrary{calc} % for manimulation of coordinates
\usetikzlibrary{shapes}
\usetikzlibrary{fit}
\usepackage{pgfplots}
\usetikzlibrary{patterns}
\usetikzlibrary{intersections,patterns,pgfplots.fillbetween}
\usetikzlibrary{decorations.pathreplacing,perspective,positioning}
\usepackage{tikz-3dplot}
\newcommand{\eg}{\textit{e.g.},~}
\newcommand{\ie}{\textit{i.e.},~}
\newcommand{\Tstrut}{\rule{0pt}{2.6ex}}       % top strut
\newcommand{\Bstrut}{\rule[-1.2ex]{0pt}{0pt}} % bottom strut

\newcommand{\Exp}{\mathbb{E}}

\newcommand{\Um}{\mathbf{U}}

\newcommand{\rhof}{\rho_{\!_f}}

\newcommand{\Prob}{\mathbb{P}}
\newcommand{\pdf}{\mathrm{p}}

\newcommand{\Thetaa}{\boldsymbol{\Theta}}

\newcommand{\sigmavm}{\boldsymbol{\sigma}_\mathrm{vm}} 
\newcommand{\sigmap}{\sigma_{\!_\mathrm{PN}}} 
\newcommand{\sigmad}{\boldsymbol{\sigma}_\mathrm{dev}}

\newcommand{\xii}{\boldsymbol\xi}

\newcommand{\rhoo}{\boldsymbol{\rho}} 
 
\newcommand{\nele}{N_\mathrm{ele}}

\newif\ifdarkmode
\darkmodefalse%\darkmodetrue
\ifdarkmode
	\definecolor{commentcolor}{rgb}{1,0,0}%{0,1,0}
	\definecolor{sourcecolor}{rgb}{1,0,0}%{1,.5,1}
	\pagecolor{black}\color{white} % dark mode
	\definecolor{newblue}{rgb}{0.7,0.7,1}
\else
	\definecolor{commentcolor}{rgb}{1,0,0}%{1,0,1}
	\definecolor{sourcecolor}{rgb}{0,.5,0}
	\definecolor{newblue}{rgb}{0,0,1}
\fi
\DeclareMathOperator*{\argmin}{arg\,min}

\title{Reliability-based Topology Optimization using Large Deviation Theory}

\author{{%\includegraphics[scale=0.06]{orcid.pdf}\hspace{1mm}
Maryam Maghazeh
}%\thanks{Use footnote for providing further
        \\
	Department of Mechanical Engineering\\
	Northern Arizona University\\
	Flagstaff, AZ 86011 \\
	\texttt{mm5335@nau.edu} \\
	\And
    Ayyappan Unnikrishna Pillai\\
	School of Infrastructure\\
	Indian Institute of Technology\\
	Bhubaneswar, India \\
	\texttt{a22ce09004@iitbbs.ac.in} \\
    \AND 
    Mohammad Masiur Rahaman\\
	School of Infrastructure\\
	Indian Institute of Technology\\
	Bhubaneswar, India \\
	\texttt{masiurr@iitbbs.ac.in} \\ 
    \AND 
    Subhayan De\thanks{Corresponding Author} \\ 
    Department of Mechanical Engineering\\
	Northern Arizona University\\
	Flagstaff, AZ 86011 \\
	\texttt{Subhayan.De@nau.edu} \\
}

\renewcommand{\shorttitle}{RBTO using Large Deviation Theory}

\begin{document}
\maketitle

\begin{abstract} 
Reliability-based topology optimization (RBTO) requires repeated estimation of small failure probabilities and their gradients, making conventional nested Monte Carlo approaches computationally prohibitive for large scale structural problems. We propose an RBTO framework that integrates large deviation theory~(LDT) with stochastic gradient descent~(SGD) to address this challenge. LDT provides closed-form exponential rate estimates of rare event probabilities, enabling accurate gradient computation without parametric assumptions on the failure density and without evaluating a full nested reliability loop at every iteration. These LDT-based gradient estimates are used directly to drive a mini batch SGD update of the design variables using only a few random samples per iteration. 
% , so that the design variables are updated using only a few random samples per iteration. 
% The probability of failure constraint is incorporated through a penalty formulation, converting the original constrained RBTO problem into an unconstrained objective that is fully compatible with the gradient based update scheme. 
The framework is validated on three benchmarks, namely, a two-dimentional (2D)  simply supported rectangular %{\color{green}(MMR: Before introducing abbreviations, full form should be written when it appears for the first time.)} 
beam, a 2D L-shaped beam, and a three-dimentional (3D) cantilever beam, under both compliance-based and stress-based failure criteria. Across these examples, the RBTO designs achieve failure probabilities lower than their robust topology optimization counterparts, demonstrating that optimizing for performance under uncertainty alone does not guarantee the satisfaction of explicit reliability constraints. Therefore, the proposed framework offers a computationally efficient route to reliable structural design under uncertainty, with direct relevance to safety-critical engineering applications. 
% Reliability-based topology optimization (RBTO) modifies topology optimization by considering uncertainties to ensure that optimized designs remain reliable under uncertain variables. In this work, we introduce a new RBTO framework that incorporates large deviation theory to calculate the probability of failure and its gradients effectively. Large deviation principles provide an efficient way to quantify how failure probabilities decay, which enables more accurate and efficient reliability estimation during the optimization process.
% Uncertainties in material properties and loading conditions are modeled using Gaussian, lognormal, and uniform distributions. Large deviation theory is used to evaluate the behavior of the probability of failure of events. Through three numerical examples, we demonstrate that the proposed approach produces designs that are both efficient and reliable. This study incorporates 2D and 3D cases and different probability of failure criteria calculation from compliance based and stress based models. It is observed that the probability of failure of robust design of 2D MBB beam, Lshaped beam and 3D cantilever beam examples are approximately $17$, $13$, $4$ times higher than that of RBTO designs, which demonstrates that robust designs are not reliable. These results suggest that the method can be implemented for applications in areas such as aerospace and civil engineering, where safety and performance under uncertainty are essential.

\end{abstract}

% keywords can be removed
\keywords{Reliability-based Topology Optimization; Large Deviation Theory; Probability of Failure; Stochastic Gradient Descent}

\section{Introduction}
% \lipsum[2]
% \lipsum[3]
Topology optimization (TO) is a computational design technique to determine the optimal distribution of material within a specified design domain, subject to applied loads, boundary conditions, and constraints, to maximize structural performance \cite{sigmund2013topology,ole-checker}. Depending on the formulation, TO can target reductions in structural compliance, weight, or improvements in performance metrics such as energy absorption, natural frequencies, thermal conductivity. Therefore, by allowing material to be freely distributed throughout the domain, TO enables the discovery of highly efficient and innovative designs that are often non-intuitive and difficult to obtain using traditional, experience-based design approaches. 
% Topology optimization (TO) is a computational method in engineering design to determine the optimal distribution of material within design domain under applied loads, boundary conditions, and constraints to maximizing the performance of the system. For example, TO can be used to minimize weight, compliance or maximizing stiffness. By topology optimization, some innovative and efficient designs can be created which might not be intuitive or achievable through traditional design methods. At first sight, It may increase the complexity of the manufacturing process; however, this concern is addressed by adding production constraints in real world problems. 
% Topology optimization process starts with domain discretizations into a finite element mesh, then solving the governing equations to find the domain response under applied loads and boundary conditions. Using these outputs as inputs for optimizer algorithm to add or remove material at each element is the next step. The final material distribution in design domain meets the performance criteria under specified constraints.

One of the primary advantages of TO is its ability to generate designs that are lighter and more efficient than those produced by conventional methods. This is especially beneficial in industries such as aerospace and automotive engineering where weight reduction is crucial for fuel efficiency and economic viability \cite{li2018multi,cavazzuti2011high,seabra2016selective,zhu2016topology}. Beyond these sectors, TO has been widely applied to the design of load-efficient structural components \cite{eschenauer2001topology,HABASHNEH2024load}, compliant mechanisms \cite{sigmund1997mechanical}, vibration-sensitive components \cite{GOLECKI2023vibration}, and heat sinks and heat exchangers \cite{dbouk2017review}. In recent years, TO has also found applications in biomedical engineering for patient-specific implants and orthopedic devices \cite{wu2021implant}. Furthermore, advances in additive manufacturing have significantly extended the practical reach of TO by enabling the fabrication of complex geometries that were infeasible \cite{gibson2015additive, wu2022topology,liu2018current}. As a result, TO has become an important design tool for multiple engineering disciplines. 

% Applications of topology optimization are widespread and include the design of mechanical components, structural parts, and microstructures. For example, in the automotive industry, it can be used to design components that are lightweight and structurally efficient to improve fuel efficiency and performance. In the field of biomedical engineering, topology optimization can help create optimized implants and prosthetics that respond to the specific needs of patients. Additionally, with the recent developments in manufacturing technologies such as additive manufacturing, the complex and innovative geometries generated by topology optimization can be more easily fabricated. This will make the optimized designs more likely to use in practical applications.

In practical engineering applications, structural performance is affected by the ubiquitous presence of uncertainty in material properties, geometric imperfections, loading conditions, and manufacturing defects. If not taken into account during the optimization process, the resulting designs can exhibit degraded performance or even failure when subjected to real-world variability \cite{pasini2019imperfect}. To address this issue, a growing body of research has been devoted to extending TO to explicitly account for uncertainty \cite{chen2010level,dunning2011introducing,maute2014topology,de2019tomaterial,KESHAVARZZADEH2017rbtomat,TOOTKABONI2012rbtomat,WU2025additive}, leading to the development of topology optimization under uncertainty (TOuU) \cite{de2019tomaterial}, where these uncertainties are often modeled using known probability distributions. 
% Moreover, in addition to the uncertainties discussed above, manufacturing defects have been incorporated into topology optimization using a deep generative framework \cite{WU2025additive}. 
% The optimization process aims to find a design that not only meets the performance criteria under nominal conditions but also remains reliable when conditions vary within acceptable limits. 

TOuU methodologies can be broadly classified into two categories, namely, robust topology optimization (RTO) (\ie without any reliability constraint) and reliability-based topology optimization (RBTO). Both approaches account for uncertainty, but they differ in how they do so. 
RTO seeks to reduce the sensitivity of the structural response to uncertain parameters, typically by minimizing a weighted combination of the mean and variance of the performance metric \citep{dunning2013robust,jansen2015robust,beyer2007robust}, without imposing an explicit constraint on the probability of failure.
RBTO, by contrast, explicitly enforces at least one constraint on the failure probability of the designed component \citep{enevo1994rbdo,royset2001rbdo,kharmanda2004reliability,lopez2012rbdo}. 
% RTO does not consider constraints on the probability of failure. On the other hand, RBTO explicitly includes at least one constraint on the failure of the designed component. In RTO, as it seeks to reduce the sensitivity of structural performance to uncertain parameters, the mean and variance of the performance are often used \cite{dunning2013robust,jansen2015robust,beyer2007robust}. 
Within RTO, different sources of uncertainty have been considered. For example, Asadpoure et al.~\cite{asadpoure2011robust} applied a perturbation method to account for
variability in material properties and structural geometry, while Guo et al.~\cite{guo2013robust} and Jansen et al.~\cite{jansen2015robust} focused on uncertainty in the geometry of the structure for robust optimization. Schevnel et al.~\cite{schevenels2011robust} represented the spatial variation of manufacturing errors using random fields. Additionally, uncertainty in the loads has also been incorporated into RTO by several authors  \cite{wang2019robust,bai2021robust,liu2018robust}. The use of surrogate models has also been explored for RTO to alleviate the computational cost associated with repeated  evaluations of the structural model \cite{KESHAVARZZADEH2017rbtomat}. 

%there has been growing interest in extending topology optimization frameworks to explicitly account for uncertainty. 

% In modern structural design, it is important to consider uncertainty while using optimization methods to make designs cost effective and reliable. Optimization approaches that account for uncertainty are generally divided into two categories, robust topology optimization and reliability based topology optimization (RBTO).
% robust topology optimization approaches rely on uncertain analyses and aim to enhance structural performance while minimizing sensitivity to variations in uncertain parameters \cite{jiansu1997rdo,wei2000rdo,messac2002rdo}.

% In contrast, RBTO approaches further enforce a constraint on satisfying a target reliability considering multiple sources of uncertainty \cite{enevo1994rbdo,royset2001rbdo,kharmanda2004reliability,lopez2012rbdo}. 
In RBTO, the target probability of failure is typically very small, which makes direct random sampling methods computationally expensive and inefficient. 
Local approximation methods, such as the reliability index approach~(RIA) \cite{haldar2000ri,hasofer1974RI,madsen2006ri,melchers2018ri} and the performance measure approach~(PMA) \cite{Tu1999pm}, are based on the first-order reliability method (FORM) and approximate the limit-state function through a first-order Taylor expansion at the most probable point (MPP) of failure in the standard normal space. 
While these approaches are computationally efficient, their reliance on local linearization can lead to significant inaccuracies in the presence of strong nonlinearities, non-convex failure domains, or multiple competing failure regions. Furthermore, when embedded within RBTO frameworks \citep{li2024level,wang2022time}, they require repeated solutions of the MPP problem at each design iteration, resulting in substantial increase in computational cost and potential convergence difficulties. 

To mitigate the cost of nested double-loop formulations, several decoupling and single-loop strategies have been proposed. Du and Chen \cite{du2004sequential} introduced Sequential Optimization and Reliability Assessment~(SORA), which separates the reliability analysis from the design update into alternating sequential sub-problems; in each cycle, the MPP identified during the reliability assessment is used to shift the deterministic constraint boundary, transforming the probabilistic constraint into an equivalent deterministic one without requiring a nested inner loop. 
Silva et al.~\cite{silva2010component} derived a single-loop RBTO formulation from Karush--Kuhn--Tucker optimality conditions, considering both component- and system-level failure probabilities. 
% To avoid this, Silva et al.~\cite{silva2010component} proposed a single-loop RBTO formulation derived from Karush--Kuhn--Tucker optimality conditions and considered both component- and system-level failure probabilities. 
Building on the SORA framework, Torii et al.~\cite{torii2016decouple} applied the decoupled sequential strategy to compliance-based robust topology optimization, and dos Santos et al.~\cite{dos2018decouple} extended it to stress-constrained RBTO using the PMA \cite{Tu1999pm}. 
% Decoupled strategies between optimization and reliability assessment were implemented in Torii et al. \cite{torii2016decouple} and dos Santos et al. \cite{dos2018decouple}. 
% To overcome these challenges, several cost reduction techniques have been developed, including multi-fidelity methods \cite{chau2019mf,gano2006mf}, response surface methods \cite{FOSCHI2002sm,agarwal2004sm}, and surrogate modeling approaches \cite{BASUDHAR2008surr,bichon2008surr,MISSOUM2007surr,sury2016surr,ZHANG2004surr,mous2019surr}. 

More recently, stochastic gradient-based formulations \cite{de2019tomaterial,de2019sgd,de2023topology} were proposed that leverage sampling techniques within the design loop for topology optimization under uncertainty, requiring only a few random samples per iteration. This approach was further extended in De et al.~\cite{de2021reliability}, who assumed a parametric model for the density of design parameters conditioned on failure and applied Bayes' theorem to derive gradient estimators for RBTO. This strategy significantly reduces the number of forward model evaluations compared to nested Monte Carlo, at the cost of assuming a specific parametric form for the conditional density and still relying on repeated limit-state evaluations. 

Despite these advances, efficiently estimating small failure probabilities and their gradients without repeated full-model evaluations and without restricting the form of the conditional failure density for nonlinear limit-state functions remains an open challenge. Large deviation theory~(LDT) \citep{amir1998ldtth,vara1984ldtthe} offers a principled analytical framework for approximating the probability of rare events through exponential rate functions, providing gradient estimates that are accurate for small failure probabilities and do not require a parametric assumption on the failure density. Its integration into RBTO has, however, not been explored. 

In this paper, we propose a stochastic gradient-based approach for RBTO that exploits LDT to efficiently compute failure probability gradients. The optimization problem, originally formulated with a constraint on the probability of failure, is transformed into an unconstrained formulation using a penalty method. Following the works of De et al.\cite{de2019sgd,de2019tomaterial} and Li et al.\cite{li2020sgd}, we employ a stochastic gradient descent~(SGD) algorithm with a small number of random samples per iteration (mini-batch) to reduce the cost. Additionally, the probability of failure is corrected at every few iterations using subset simulation \cite{AU2001loc} to ensure accuracy without prohibitive computational cost. Gradient estimation of the failure probability is performed via LDT, providing an effective analytical route to account for rare-event probabilities in reliability-based design without a full nested Monte Carlo loop. The method is demonstrated on two-dimensional and three-dimensional examples with different constraints, showing a practical approach for real engineering problems.

% originally formulated with a constraint on the probability of failure then is transformed into an unconstrained formulation using a penalty method.
% Following the works of De et al.\cite{de2019sgd,de2019tomaterial} and Li et al.\cite{li2020sgd}, we employed a stochastic gradient descent (SGD) algorithm with a small number of random samples per iteration called mini batch to reduce the computational cost of failure probability estimation. Specifically, the probability of failure is evaluated every few iterations using subset simulation to ensure accuracy without high computational costs.
% Since the proposed framework relies on gradient-based optimization, we incorporate the concept of Large Deviation Theory (LDT) \cite{amir1998ldtth,vara1984ldtthe} to compute the gradient of the probability of failure efficiently. We note that the development of efficient techniques for estimating probability of failure and corresponding gradients, particularly for small failure probabilities, such as large deviation theory is an interesting topic for researchers. The integration of LDT into RBTO represents a novel contribution, providing an effective way to account for probabilities of events in reliability based design. The method successfully incorporated in 2D and 3D examples with different constraints showing a practical approach for real engineering problems.

The remainder of the paper is organized as follows. Section \ref{sec:background} presents the necessary background on density-based topology optimization using the solid isotropic material with penalization (SIMP) 
% {\color{green}(MMR: Before introducing abbreviations, full form should be written when it appears for the first time.)} 
approach, stress-based formulations, and the fundamentals of large deviation theory. Section \ref{sec:method} introduces the proposed LDT-based RBTO formulation, the stochastic gradient descent update scheme, and the overall algorithmic framework. Thereafter, Section \ref{sec:examples} provides three numerical examples, including two-dimensional simply supported rectangular and L-shaped beam problems and a three-dimensional cantilever beam, demonstrating the effectiveness of the method for both compliance-based and stress-based failure criteria. Finally, Section \ref{sec:conclusions} summarizes the main findings and discusses the computational advantages and practical implications of the proposed approach for large-scale reliability-based structural design.

\section{Background} \label{sec:background}

% This section provides the background on topology optimization. Subsequently, large deviation theory used in the proposed method is briefly discussed. %With this background, the proposed use of 
Addressing the presence of uncertainties in loading, material properties, and manufacturing processes during topology optimization requires a probabilistic treatment of rare but critical failure events. Therefore, this section reviews the theoretical foundations of both topology optimization and large deviation theory, which together form the basis of the proposed method.

\subsection{Topology Optimization}

We employ the SIMP 
% {\color{green}(MMR: This abbreviation appears in the previous section.)} 
method \cite{bendose1989simp,ZHOU1991simp} for material interpolation within the topology optimization framework. 
The SIMP method employs a power-law scheme to relate the effective Young's modulus  $E(\rho)$ to the local density variable $\rho\in[0,1]$, where $\rho = 1$ denotes fully solid material and $\rho = 0$ denotes void. The interpolation is expressed as  
% The SIMP method uses a power-law interpolation scheme to model the dependence of material stiffness on the density variable $\rho\in[0,1]$ with $\rho=1$ representing a solid region, and $\rho=0$ representing a void region. The effective Young's modulus $E(\rho)$ is interpolated as:\ 
\begin{equation}
        E(\rho)=(\eta_{e} + (1-\eta_{e})\rho^q)E_0, 
\end{equation} 
% where  $q$ is a penalization factor, generally recommended as $q\geq3$, to suppress intermediate densities and obtain a clear optimized design; $E_0$ is the modulus of elasticity of the solid material; and $\eta_{e}=10^{-15}$ is a small number introduced to prevent singularity in void regions. 
where $q\geq3$ 
% {\color{green}(MMR: How is it meaningful?)} 
is the penalization factor, which suppresses intermediate densities to yield a well-defined solid-void design, $E_0$ is the Young's modulus of the solid material, and $\eta_{e} = 10^{-15}$ is a small positive constant introduced to prevent stiffness matrix singularity in void regions. 
% The material densities are defined in terms of the design parameters $\rho$ using a Helmholtz-type density filter subject to homogeneous Neumann boundary conditions, which enforce zero flux across the boundary. \cite{ole2011helm}. This PDE-solved filtering step helps regularize the design and avoid numerical instabilities such as checkerboarding \cite{ole-checker}. 
To ensure mesh-independence and suppress numerical instabilities such as checkerboarding \cite{ole-checker}, the design variables are regularized through a Helmholtz-type density filter with homogeneous Neumann boundary conditions \cite{ole2011helm} as 
\begin{equation}
    \begin{split}
    -r_f^2\nabla^2\rhof + \rhof = \rho \quad \text{in} \: \Omega, \\
    \nabla\rhof \cdot\mathbf{n} = 0 \quad \text{on} \: \partial\Omega,
    \end{split}
\end{equation}
% {\color{green}(MMR: Please check the above equation.)}
% where $\rho_f$ is filtered density and $r_f$ is filter radius which set as $r_f = 3h_c/ 2\sqrt{3}$ \cite{ole2011helm}. $h_c$ representing the characteristic element size. 
where $\rhof$ is the filtered density field, $r_f = 3h_c / 2\sqrt{3}$ is the filter radius, and $h_c$ is the characteristic element size. The Neumann boundary condition enforces zero flux across the domain boundary, ensuring that the filter does not introduce artificial density gradients near the edges.
% To promote a discrete, we additionally apply a projection scheme to the filtered design variables, enabling clearer separation between solid and void regions \cite{ole2011thresh}. 
To recover a near-binary material distribution from the filtered field, a smooth projection scheme is subsequently applied \cite{ole2011thresh} as
\begin{equation}
    \rho_t = \frac{\tanh(\beta\eta)+\tanh(\beta(\rhof-\eta))}{\tanh(\beta\eta)+\tanh(\beta(1-\eta))}
\end{equation}
% where $\rho_t$ are the threshold density $\beta$ corresponds to the projection strength and $\eta$ is the projection threshold parameter. We used $\eta = 0.5$ for all the examples. 
where $\rho_t$ is the projected density, $\beta$ controls the sharpness of the projection, and $\eta$ is the threshold parameter defining the transition point between solid and void. A projection threshold of $\eta = 0.5$ is adopted herein, following standard practice in the literature. 

\subsection{Stress-based Topology Optimization} 
% Topology optimization involving compliance and volume as quantity of interest is well established with straightforward gradient calculations. However, stress-based topology 
Stress exceeding a prescribed safe level is one of the primary structural failure criteria in engineering design. Failure criteria based on stress constraints can be included in the RBTO formulation, as is used in Example II in Section \ref{sec:examples}, to ensure structural integrity under uncertainty. 
% In the context of RBTO, this motivates the imposition of stress constraints to ensure structural integrity under uncertainty. 
However, because stress is a local quantity, enforcing pointwise stress constraints \cite{duysinx1998stress,bruggi2012stress} generates a large number of local constraints, rendering large-scale problems computationally intractable. 
% Duysinx et al.~\cite{duysinx1998stress} and Bruggi et al.~\cite{bruggi2012stress} demonstrated the feasibility of local stress constraints in small-scale problems. 
To address this, local stress values are commonly aggregated into a single global measure, with the $p$-norm being a widely adopted approach due to its smoothness and  differentiability for gradient-based optimization~\cite{yang1996stress,picelli2018stress,le2010stress, GU2024pnorm}. 
% {\color{green}(SD: The second reference is a random reference. There are many other famous references on p-norm. Replace that. You can also use \href{https://scholar.google.com/scholar_labs/search?hl=en}{Google Scholar Lab})}{\color{green}(MM: I added a paper with more citations. the last paper was one of the latest papers but less citation)}
% One of the structure failure criteria can be stress exceeding a safe stress level, reliability based topology optimization based on the stress criteria is discussed in this section.
% As stress is a local quantity, stress-based RBTO requires many local constraints which is computationally expensive and makes large-scale problems impractical. Duysinx et. al.~\cite{duysinx1998stress} and Bruggi et. al.~\cite{bruggi2012stress} applied local stress constraints in smaller studies successfully 
% To address this issue, local stress criteria are often aggregated into global stress measurements. $P$-norm stress is a common aggregation method to approximate the maximum stress \cite{le2010stress,GU2024pnorm} while remaining smooth and differentiable for gradient-based optimization. 

Since this study considers linear elasticity for homogeneous isotropic material, the von Mises stress is adopted as the local stress measure and defined as 
\begin{equation}
    \sigmavm = \sqrt{\frac{3}{2}\sigmad:\sigmad},
\end{equation}
where $\sigmad$ is deviatoric part of the Cauchy stress tensor. 
The local von Mises stress values are then aggregated into a global $p$-norm stress measure, defined as
% The global $p$-norm stress is defined as
\begin{equation}\label{eq:pnorm}
    \sigmap = \left(\sum_{n=1}^{\nele}\sigmavm^p\right)^{1/p},
\end{equation} 
where $p$ is the norm aggregation parameter, and $\nele$ is the number of elements in the structural model. Higher values of $p$ yield a closer approximation to the true maximum stress over the domain, but can cause numerical ill-conditioning %\cite{ole2021pnorm}. 
% is the penalization factor for $p$ norm function. The higher values of $p$ shows a better stress distribution but makes the problem unstable 
\cite{ole2021pnorm,de2015pvalue,le2010stress}. 
% {\color{green} (SD: Again, the first one here is a random reference. We are not even doing binary variables.)} 
In the limit $p \to \infty$, converges to the global maximum von Mises stress. A value of $p = 30$ is adopted herein, which is consistent with recommendations in the literature \cite{lian2017pvalue}. %\cite{zhou2017pvalue,lian2017pvalue}.
The stress constraint used in the RBTO formulation is therefore expressed as 
% The von Mises stress is used as this study considers linear elasticity theory for homogeneous material.
% The von Mises stress defined as 
% \begin{equation}
%     \sigmavm = \sqrt{\frac{3}{2}\sigmad:\sigmad}
% \end{equation}
% where $\sigmad$ is deviatoric part of the stress. Therefore, the criteria used in stress based RBTO defined as 
\begin{equation}
    \sigmap \leq \sigma_\mathrm{yield}
\end{equation} 
where $\sigma_\mathrm{yield}$ is the prescribed yield stress of the material.

Since we employ a gradient-based optimization approach, the required sensitivity of $\sigmap$ with respect to the displacement field $\Um$ can be estimated by applying the chain rule as
% the gradients of this function must be computed. The gradients of the $p$-norm stress with respect to displacement $\Um$ can be calculated by the chain rule as follows:
\begin{equation}
    \frac{\partial \sigmap}{\partial \Um} = \frac{\partial \sigmap}{\partial \sigmavm}\frac{\partial \sigmavm}{\partial \sigmad}\frac{\partial \sigmad}{\partial \Um}, 
\end{equation}
% The  von Mises stress at the center of the $i$th element defined as 

% ${\sigma}_{\mathrm{vm},i}=\left(\sigma_{x,i}^2+\sigma_{y,i}^2-\sigma_{x,i}\sigma_{y,i}+3\tau_{xy,i}^2\right)^{1/2}$ where, $\sigma_{x,i}$ and $\sigma_{y,i}$ are axial stresses in $x$- and $y$-direction, respectively, $\tau_{xy,i}$ is the shear stress, and $\nele$ is the total number of elements. 
% \begin{equation}
%     \sigma_{vm,i} = \left(\sigma_{ix}^2+\sigma_{iy}^2-\sigma_{ix}\sigma_{iy}+3\tau_{xy}^2\right)^{1/2} 
% \end{equation}
% we penalized each $\sigma$ by a penalization factor $\eta$
% \begin{equation}
%     \hat{\sigma_i}(\rho) = \eta(\rho)\sigma
% \end{equation}
% We used the penalization scheme as $\eta(\rho) = \rho^q$ 
% \subsubsection{Sensitivity Analysis} 
% Gradients calculation as a important step for topology optimization is explored in this section.
where the first two terms can be expressed as 
\begin{equation}
\begin{split}
    &\frac{\partial \sigmap}{\partial \sigmavm} = \left(\left(\sum_{n=1}^{\nele}\sigmavm^p\right)^{1/p -1}\right)\sigmavm^{p-1};\\
    &\frac{\partial \sigmavm}{\partial \sigmad} = \frac{3}{2\sigmavm}(\sigmad).\\
    % &\frac{\partial \sgvm}{\partial \sigma_{ix}} = \frac{1}{2\sgvm}(2\sigma_{ix} - \sigma_{iy})\\
    % &\frac{\partial \sgvm}{\partial \sigma_{iy}} = \frac{1}{2\sgvm}(2\sigma_{iy} - \sigma_{ix})\\
    % &\frac{\partial \sgvm}{\partial \tau_{xy,i}} = \frac{3\tau_{xy,i}}{\sgvm}\\
\end{split}
\end{equation}
% In Section \ref{sec:examples}, one of the numerical example uses a stress based failure criteria. 
% The adjoint vector is then calculated based on this result. 

\subsection{Large Deviation Theory} \label{sec:ldt} 
LDT 
% {\color{green} (MMR: This abbreviation appears in Section 1 also.)}{\color{green} (MM: so i should just use LDT?)}
is a branch of probability theory that studies the asymptotic behavior of rare events. It provides tools to quantify the probability of deviations of a random variable from its typical or expected value. The framework is particularly useful for understanding the probability of events that occur with exponentially small likelihood, such as the rare failures in engineering structures. Using LDT, the log of probability of failure constraint can be written as 
\begin{equation}
    \begin{split}
        &\log\Prob(g(\rhoo;\xii)\geq z) \approx -I(\xii^*(\rhoo,z))\qquad\text{as } z\rightarrow\infty, 
    \end{split}
\end{equation} 
where $I(\cdot)$ is the rate function, $\xii$ is the vector of uncertain variables, 
% {\color{green}(I added these)}
$g(\rhoo;\xii)$ is the performance function, $z$ is the failure limit, $\xii^*(\rhoo,z)$ is the MPP 
% {\color{green} (MMR: Earlier the same abbreviation appears for most probable failure point.)} 
or \textit{instanton} given by 
\begin{equation}
    \begin{split}
        &\xii^*(\rhoo,z):=\argmin_{\xii\in\Xi} \{I(\xii): g(\rhoo;\xii)\geq z\}. 
    \end{split}
\end{equation} 
\begin{figure}[!htb]
    \centering
    \includegraphics[width=0.5\linewidth]{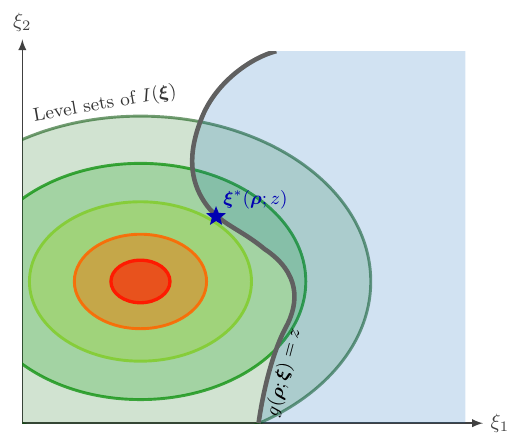}

    \caption{Most probable point $\xii^*(\rhoo;z)$ from large deviation theory.}
    \label{fig:ldt_schem}
\end{figure}
The rate function encapsulates the essential information about how rare events occur and how their probabilities decay. The importance of the rate function is that it simplifies the analysis of probabilities for rare events. Instead of calculating probabilities directly (which is often intractable for small probabilities), the rate function gives an asymptotic estimate. This property makes LDT a powerful tool for understanding systems with low-probability, high-impact events. 
% Table \ref{tab:rate_fun} shows the rate functions for a few typical probability distributions. 
% {\color{red}(SD: Need a few sentences connecting $\xii^*$ with FORM and $t^*$.)}{\color{red}(mm:t* is not defined here)}
While closed form expressions are available for the rate function for typical probability distributions, \eg Gaussian, exponential, and so on, for a general continuous random variable, the rate function can be derived using the Legendre-Fenchel transform as \cite{fenchel1949conjugate} 
\begin{equation} \label{eq:rate_fun_sup}
    I(\xii) = \sup_{t\in\mathbb{R}} [t\xii-\Lambda(t)], 
\end{equation} 
where $\Lambda(t)=\log\Exp[\exp(t\xii)]$ is the cumulant generating function (CGF) provided the moment-generating function exists. This optimization problem can be interpreted as the dual of the MPP formulation, with the optimization carried out in the space of exponential tilting parameters induced by the CGF rather than in the original random variable space. 
% The CGF encapsulates the exponential moments of the distribution and plays a central role in deriving rate functions via convex duality. 
The supremum in \eqref{eq:rate_fun_sup} is attained at $t^*$, which satisfies the stationarity condition $\xii = \Lambda'(t^*)$. Consequently, the rate function can be written as $I(\xii) = t^* \xii - \Lambda(t^*)$.
Note that, $t^*$ is obtained through the CGF, which enables a global characterization of the rare event and avoids reliance on local approximations of the performance function. As a result, the LDT-based formulation exhibits improved robustness in the presence of strong nonlinearity. %Moreover, invoking the contraction principle,

\section{Methodology} \label{sec:method}
In this section, we define the reliability-based topology optimization problem followed by the proposed use of large deviation theory along with a stochastic gradient descent approach for the solution of this problem. The proposed method is then summarized in Algorithm \ref{alg:ldt-rbto}. 
% \subsection{Topology}

\subsection{Reliability-based Topology Optimization (RBTO)} 

The reliability-based topology optimization problem used herein is defined as 
\begin{equation} \label{eq:problem}
\begin{split}
    &\min_{\rhoo\in\Thetaa} \omega_c\Exp\left[C(\rhoo;\xii)\right] + \omega_v V(\rhoo) \\
    &\text{subject to } K(\rhoo;\xii)U(\rhoo;\xii) = F(\xii) \\ 
    &\qquad\qquad P_f:=\Prob(g(\rhoo;\xii)\geq z ) \leq p_a, \\  
\end{split} 
\end{equation} 
where $C(\rhoo;\xii)$ is the compliance of the design, $V(\rhoo)$ is the total volume of the design, $\omega_c$ and $\omega_v$ are the weights for the compliance and volume, respectively, $K(\rhoo;\xii)$ denotes the stiffness of the structure, $U(\rhoo;\xii)$ denotes the displacement of the structure, and $F(\xii)$ is the uncertain force. The failure probability $P_f$ is defined as the probability of a given 
% {\color{green} (MMR: What do you mean by some? Be specific.)} 
performance metric $g(\rhoo;\xii)$ exceeding a prescribed threshold $z$ that is restricted to below an acceptable probability of $p_a$. 
Using a penalty approach, the above constrained optimization can be converted to an unconstrained optimization problem as 
\begin{equation} \label{eq:opt_uncon}
    \min_{\rhoo\in\Thetaa}~~J(\rhoo) :=  \omega_c\Exp\left[C(\rhoo;\xii)\right] + \omega_v V(\rhoo) + \frac{\kappa_f}{2} \left[\left(\ln P_f -\ln p_a \right)^+\right]^2, 
\end{equation} 
where $\kappa_f$ is a penalty parameter and the constraint violation is denoted by $(\cdot)^+$ with $(\cdot)^+=0$ for $(\cdot)<0$. 
% {\color{red} (MMR: Is it consistent with the constraint in Eq. (12)?.explain the reason for ln)}. 
% {\color{blue}(
Note that as the changes in $P_f$ can be in orders of magnitude, the natural log of $P_f$ is used in the unconstrained objective \eqref{eq:opt_uncon}. 
% )}. 
The gradient of the objective can be computed as 
\begin{equation} \label{eq:grad_obj}
    \nabla_{\rhoo} J(\rhoo)= \omega_c\Exp[\nabla_{\rhoo} C(\rhoo;\xii)] +  \omega_v\nabla_{\rhoo} V(\rhoo) + \kappa_f \left[(\ln P_f-\ln p_a)^+\right] \left(\nabla_{\rhoo}\ln P_f\right), 
\end{equation} 
where the gradients of the compliance and volume can be obtained in a straightforward manner \cite{ole-checker}. 
% following Andreassen study \cite{Andreassen2011gradcv}. 
However, \eqref{eq:grad_obj} requires the gradient of the log of probability of failure and poses a challenge. The large deviation theory is used to efficiently estimate this gradient as discussed in the next subsection. 

Once the gradient of the objective function, $\nabla_{\rhoo} J(\rhoo)$ is estimated, a stochastic gradient descent method \cite{de2019sgd,de2021reliability} is used for a scalable solution of \eqref{eq:opt_uncon}, which updates the density variable at $k$th iteration as 
\begin{equation} \label{eq:sgd_to}
    \rhoo_{k+1} = \rhoo_k-\eta \nabla_{\rhoo}J(\rhoo_k), 
\end{equation} 
where $\eta$ is a step size or learning rate parameter. Herein, a mini-batch stochastic gradient descent approach is employed using a handful of samples per iteration instead of using a large number of random samples to approximate the gradients in \eqref{eq:sgd_to}.

% \subsection{Comparison with First Order Reliability Method (FORM)} 

\subsection{Proposed Use of Large Deviation Theory}\label{subsec:Proposed Use of Large Deviation Theory}
% \subsubsection{Approach I: LDT-SORA} 

% In the first approach, we cast the RBTO problem using the SORA method by decoupling the nested optimization problem that involves estimating $\xii^*$. 
% \begin{equation}
% \begin{split}
%     & \textbf{Deterministic Optimization:}\\ 
%     &\min_{\rhoo\in\Thetaa} V(\rhoo) \\
%     &\text{subject to } K(\rhoo;\xii^*)U(\rhoo;\xii^*) = F(\xii^*) \\ 
%     &\phantom{\text{subject to } } g(\rhoo;\xii^*)\geq z \\ 
%     & \textbf{Rate Function Estimation:}\\ 
%     &\xii^*=\argmin_{\xii\in\Xi} I(\xii) \\ 
%     &\text{subject to } g(\rhoo;\xii)\geq z \\ 
% \end{split} 
% \end{equation} 

% \subsubsection{Approach II: Estimated CGF-based RBTO} 

In the proposed approach, we perform sample-based estimates of the CGF of the performance function, \ie $g(\rhoo;\xii)$, such as the compliance or maximum von Mises stress. 
The empirical CGF using random samples of $\{g(\rhoo;\xii_i)\}_{i=1}^n$ can be estimated as  
\begin{equation}
    \widehat{\Lambda}_g(t; \rhoo) = \log\left( \frac{1}{n} \sum_{i=1}^n \exp\left({t g(\rhoo;\xii_i)}\right) \right). 
\end{equation} 
While this is a biased estimator of the true CGF because of the log outside the expectation, but it is consistent as $n \to \infty$.
The solution of \eqref{eq:rate_fun_sup} is given by 
\begin{equation} \label{eq:tilt}
    \frac{\mathrm{d}}{\mathrm{d}t}\left(\Lambda_g(t; \rhoo)\right) = \Exp_\mathrm{tilted}\left[g(\rhoo;\xii)\right] = z, 
\end{equation} 
where the \textit{tilted expectation} is with respect to a exponentially tilted probability density $\pdf_t(x)=\pdf(x)\exp(tx)/\Exp[\exp(tx)]$, 
% {\color{green}(shall we write it in inline format??)},
and it biases the samples toward the rare event region. Here, $\pdf(x)$ is the probability density function of $g(\rhoo;\xii)$, and $t$ is the tilting parameter. 
Using samples of $g(\rhoo;\xii)$, \eqref{eq:tilt} can be solved for $t^*$ as the solution of \eqref{eq:rate_fun_sup} using the approximation of the tilted expectation,  
\begin{equation}
    \Exp_\mathrm{tilted}\left[g(\rhoo;\xii)\right] \approx \sum_{i=1}^n w_i g(\rhoo;\xii_i), 
\end{equation} 
where 
% {\color{green}(shall we write it in inline format??)} 
$w_i=\exp({t g(\rhoo;\xii_i)})/\sum_{j=1}^n\exp({t g(\rhoo;\xii_j)})$. 
Herein, Newton-Raphson method \cite{ypma1995historical} is used to find $t^*$. With $t^*$, the rate function can be given by 
\begin{equation} 
    I(z;\rhoo) \approx t^*z-\widehat{\Lambda}_g(t^*;\rhoo). 
\end{equation} 
Note that while the rate function was introduced in a generic form $I(\xii)$ in Section \ref{sec:ldt}, the present setting involves a design-dependent performance function $ g(\rhoo;\xii)$. Hence, the rate function is expressed as $I(z;\rhoo)$ to emphasize its parametric dependence on $\rhoo$. 
With this estimate of the rate function, the reliability constraint can be differentiated to give 
\begin{equation}
\begin{split}
    \nabla_{\rhoo} \ln P_f(\rhoo) &= - \nabla_{\rhoo} I(z;\rhoo)\\ 
    &= \nabla_{\rhoo}\Lambda_g(t^*;\rhoo)\\
    &= \nabla_{\rhoo}\log\Exp[\exp(t^*g(\rhoo;\xii))]\\
    &= \frac{1}{\Exp[\exp(t^*g(\rhoo;\xii))]}\Exp\left[t^*\frac{\mathrm{d}g(\rhoo;\xii)}{\mathrm{d}\rhoo}\exp(t^*g(\rhoo;\xii))\right]\\
    &= t^* \Exp_\mathrm{tilted}\left[\frac{\mathrm{d}g(\rhoo;\xii)}{\mathrm{d}\rhoo}\right],\\ 
\end{split}
     % = t^* \Exp_\mathrm{tilted}\left[\frac{\mathrm{d}X(\rhoo)}{\mathrm{d}\rhoo}\right],  
\end{equation} 
where the tilted expectation $\Exp_\mathrm{tilted}$ is defined as before. %the last expression follows from Appendix A. 
Therefore, using realizations of $g(\rhoo;\xii)$, the gradient can be approximated as 
% \begin{equation}
% \begin{split}
%     \nabla_{\rhoo} I(z;\rhoo)
%     &= -\nabla_{\rhoo}\Lambda(t^*;\rhoo)\\
%     &= -\nabla_{\rhoo}\log\Exp[\exp(t^*X(\rhoo))]\\
%     &= -\frac{1}{\Exp[\exp(t^*X(\rhoo))]}\Exp\left[t^*\frac{\mathrm{d}X(\rhoo)}{\mathrm{d}\rhoo}\exp(t^*X(\rhoo))\right]\\
%     &= -t^* \Exp_\mathrm{tilted}\left[\frac{\mathrm{d}X(\rhoo)}{\mathrm{d}\rhoo}\right] 
% \end{split}
% \end{equation} 
\begin{equation}
    \nabla_{\rhoo} \ln P_f(\rhoo) \approx t^* \sum_{i=1}^N w_i \frac{\mathrm{d}g(\rhoo;\xii_i)}{\mathrm{d}\rhoo}. 
\end{equation} 
These steps are summarized in Algorithm~\ref{alg:ldt-rbto}.

\paragraph{Remarks} 
\begin{itemize} 
    \item Due to the small number of samples drawn at each iteration, it is possible that no failure samples are observed at certain iterations. In such cases, the tilting parameter $t^*$ associated with the failure event cannot be identified, which can lead to a breakdown in the gradient computation of the probability of failure. To address this, samples are accumulated over every $n_s$ iterations to avoid any large memory requirements and reused in iterations where no failure samples are observed. Furthermore, we reuse the latest gradient of the probability of failure, if needed, to avoid any numerical issues. This strategy ensures the robustness of the gradient computation while preserving the computational efficiency of the stochastic optimization scheme. 
    % As mentioned before, only a small number of samples are used at each iteration. Therefor, it is possible that no failed samples are detected at certain iterations. This situation can prevent the identification of the tilting parameter $t^*$ leading to failure in computing the gradient of the probability of failure. To address this issue, we store samples over every $n_{th}$ iteration and reuse these samples when needed. This method ensures the robustness of the gradient computation while maintaining the computational efficiency of the stochastic optimization. 
    % {\color{green}(SD: I don't see any $n_s$ mentioned in the Algorithm.)} 
    \item To correct the estimate of the probability of failure, we use an efficient algorithm, such as subset simulation \cite{AU2001loc} used herein (see Appendix A), at every $m$ iterations to estimate $P_f$. The update frequency interval $m$ is selected to balance the additional computational cost involved with it and the accuracy in $P_f$. In the numerical examples, we used $m=20$, which proved to be sufficient. 
\item 
A further challenge in the unconstrained optimization formulation in \eqref{eq:problem} is the selection of appropriate weights (\ie $\omega_c$ and $\omega_v$) for the competing objective terms. If not selected judiciously, it can cause one term to dominate the objective, effectively suppressing the influence of the remaining terms on the optimized design. For example, an excessively large weight on the volume term drives the optimizer to remove material aggressively. This will prevent the reliability constraint from being enforced, as a result, yielding infeasible designs. In this paper, for each of the numerical examples, we chose these values of the weights to produce meaningful designs. However, we kept these values the same when comparing different methods. 
% Another challenge in unconstrained optimization is selecting appropriate weights for the different terms in the objective function. To ensure that each objective term influences the design, some weight combinations were tested. For example, when the volume term dominated the objective, the optimizer tended to remove excessive material, which prevented the reliability term from being properly enforced. 
\end{itemize}

\begin{algorithm}[!htb] 
\caption{Large deviation theory based reliability based topology optimization (LDT-RBTO)} \label{alg:ldt-rbto}
\begin{algorithmic}[1] 
\State \textbf{Input:} Initial topology $\rhoo_0$, %volume fraction constraint $\bar{\rho}$, 
% threshold $z$, number of samples $N$, maximum iterations $K$, 
number of samples per iteration $n$, sample accumulation frequency $n_s$, 
target failure probability $p_a$, $P_f$ update frequency $m$, penalty parameter $\kappa_f$, learning rate $\eta$, maximum number of iterations $N_\mathrm{it}$ 
% \For{$k = 0,\dots$}
\For{$k = 1$ to $N_\mathrm{it}$} 
    \State Generate $n$, a small number of \textit{i.i.d.} samples $\{\xii_i\}_{i=1}^n$ of the uncertain parameters from the pdf $\pdf(\xii)$
    \State Reset failure samples history at every $n_s$ iterations. % to address the issue related to the first remark in \ref{subsec:Proposed Use of Large Deviation Theory} %{\color{red}(SD: Need to have a discussion.)}{\color{red}(mm:I added the section which is discussed about this)}
    % \State Append the current samples in the samples histories.
    \For{$i = 1$ to $n$}
        \State Solve $K(\rhoo_k;\xii_i)U(\rhoo_k;\xii_i) = F(\xii_i)$ to obtain displacement $U(\rhoo_k, \xii_i)$
        \State Evaluate performance function $ g_i := g(\rhoo_k;\xii_i)$ 
        \State If $g_i\geq z$, append the current sample to the failure samples history.
    \EndFor
    \State Compute empirical CGF
    \[
        \widehat{\Lambda}_g(t; \rhoo_k) = \log\left( \frac{1}{n} \sum_{i=1}^n e^{t g_i} \right)
    \]
    \State Find saddle point $t^*$ such that 
    \[
        z = \frac{ \sum_{i=1}^n g_i e^{t^* g_i} }{ \sum_{i=1}^n e^{t^* g_i} }
    \]
    \State If $t^*$ cannot be identified using current samples, use the accumulated samples to find the $t^*$ (see the first remark in Section~\ref{subsec:Proposed Use of Large Deviation Theory}); %{\color{blue}(If $t^*$ cannot be identified using current samples, use the accumulated samples to find the $t^*$; otherwise, use the latest reliability gradients calculated using $t^*$.)}

    % \State Compute rate function:
    % \[
    %     I(z; \rhoo_k) = t^* z - \tilde{\Lambda}_X(t^*; \rhoo_k)
    % \]
    % \State Estimate probability of failure:
    % \[
    %     P_f(\rhoo_k) \approx \exp(-I(z; \rho^k))
    % \]
    \State Compute tilted weights
    \[
        w_i = \frac{e^{t^* g_i}}{\sum_{j=1}^n e^{t^* g_j}}
    \]
    \State Compute sensitivity
    \[
        \nabla_{\rhoo} \ln P_f(\rhoo_k)=-\nabla_{\rhoo} I(z; \rhoo_k) = t^* \sum_{i=1}^n w_i \cdot \nabla_{\rhoo} g_i
    \]

% ({\color{green} MMR: In some expression $\rho_k$ is bold (for instance in, $P_f(\rhoo_k)$) but in some expression it is non-bold (for instance in, $I(z; \rho_k)$).}

    \If{$k/m$ is an integer}
        \State Probability of failure $P_f$ calculation using efficient methods such as subset simulation (see Algorithm~\ref{alg:subset}).
    \EndIf

    % \State Update topology using gradient-based optimizer (e.g., OC or MMA):
    % \[
    %     \rho^{k+1} = \text{Update}(\rho^k, \nabla_\rho I, \bar{\rho})
    % \]
    \State Estimate the gradient of the unconstrained objective 
    \[
    \nabla_{\rhoo} J(\rhoo_k)= \omega_c \left(\sum_{i=1}^n\nabla_{\rhoo} C(\rhoo_k;\xii_i)\right) +  \omega_v\nabla_{\rhoo} V(\rhoo_k) + \kappa_f \left[(\ln P_f(\rhoo_k)-\ln p_a)^+\right] \left(\nabla_{\rhoo}\ln P_f(\rhoo_k)\right)
    \] 
    \State Update the density variables 
    \[
    \rhoo_{k+1} = \rhoo_k-\eta \nabla_{\rhoo}J(\rhoo_k) 
    \]
\EndFor
\State \textbf{Output:} Optimized topology $\rhoo_k$
\end{algorithmic}
\end{algorithm} 

\subsection{Computational Benefits} 
The primary computational cost of RBTO is from the evaluation of the probability of failure and its gradients. % represents one of the most computationally demanding components of the RBTO framework. 
To reduce this cost, three complementary strategies are employed herein. 
\begin{itemize}
    \item First, the gradient of the probability of failure is estimated using efficient sampling methods, namely, the large deviation theory, in place of direct Monte Carlo simulation, removing the requirement of a prohibitively large number of samples for accurate rare event estimation. 
    \item Second, any error in the failure probability estimation is corrected only at every $m$-th iteration by using subset simulation method \cite{AU2001loc}. A value of $m = 20$ is adopted herein, providing a practical balance between computational efficiency and solution accuracy. 
    \item Third, an SGD scheme is employed, wherein a small number of random samples per iteration is used for gradient estimation, substantially reducing the per iteration computational cost. 
\end{itemize}

% Calculation of the probability of failure and its gradients is one of the computationally costly steps in RBTO cases. To address this, we combine efficient sampling strategy as Large Deviation Theory and subset simulation instead of Monte Carlo simulation which needs a large number of samples to calculate the probability of failure. To further reduce computational expense, the failure probability is evaluated at every $m$ iterations. In our examples $m=20$ to have a balance between accuracy and efficiency. Moreover, using the stochastic gradient descent scheme, with small number of random samples per iteration for gradient calculation helps us to decrease the computational costs. 

\section{Numerical Illustrations}\label{sec:examples} 
Three numerical examples are used in this section to illustrate the efficacy of the proposed approach. The first two examples focus on two-dimensional structural design problems involving the MBB and L-shaped beams. The third example considers a three-dimensional cantilever beam. 
% As mentioned in Methodology section, the optimization objective is formulated as a weighted combination of structural compliance and material volume, and a penalty term associated with the probability of failure. Uncertainty is incorporated in both external loading and material properties to consider real-world scenarios. 
% In all three examples, consistent units are used. 
Unless otherwise mentioned, we use consistent units in the numerical examples. 
The densities are initialized to 1.0 for all elements, so that the reliability constraint is satisfied initially. The optimization is performed using the SGD approach with a learning rate of $\eta=0.075$ and a mini-batch size of $n=10$. These values provide convergence of the optimization within reasonable number of iterations. %Each iteration uses a mini-batch of $n=10$ random samples, which makes the proposed model efficient.{\color{green}(mm: added for comment in conclusion)} 
% To estimate the failure probability during the optimization process, subset simulation is employed followed the method of De et.al\cite{de2019sgd}. This framework efficiently captures rare-event probabilities and ensures that the reliability penalty is accurately incorporated into the objective function throughout the optimization. The probability of failure is estimated using the subset simulation at every $m=20$ iterations. 
In the subset simulation, which is used at every $m=20$ iterations to estimate $P_f$, the conditional probability level is set to $p_0 = 0.1$, and estimated with $N_s = 200$ samples. % are used per subset level. 
The proposed optimization approach is carried out until no significant changes in either the optimized design or the probability of failure are observed, indicating convergence. 
% conditional probability is set to $p_0=0.1$, and the number of samples per level is $N=200$. 

\subsection{Example I: Design of a Two-dimensional Rectangular Beam} 
\begin{figure}[!htb]
    \centering
    % 	\centering
    % \begin{subfigure}[t]{\textwidth}
    % \centering
    % \includegraphics[scale=1]{figs/Fig4a.pdf}
	\begin{tikzpicture}[scale=0.75,every node/.style={minimum size=1cm},on grid]
	\draw [ultra thick,fill=gray!40,draw=none] (0,-3.5) rectangle (10,0);
	\draw [ultra thick] (0,0) -- (10,0);
	\draw [ultra thick] (10,0) -- (10,-3.5);
	\draw [ultra thick] (10,-3.5) -- (0,-3.5);
	\draw [ultra thick] (0,-3.5) -- (0,0);
	
	% Rollers
	\draw[ultra thick,fill=gray!30] (9.85,-4.15) circle (0.125);
	\draw[ultra thick,fill=gray!30] (10.15,-4.15) circle (0.125);
	\draw[fill=gray!30]    (10,-3.5) -- ++(0.3,-0.5) -- ++(-0.6,0) -- ++(0.3,0.5);
	% Triangle
	\draw [ultra thick] (10,-3.5) -- (10.3,-4);
	\draw [ultra thick] (10.3,-4) -- (9.7,-4);
	\draw [ultra thick] (9.7,-4) -- (10,-3.5);
	
% 	% Rollers
% 	\draw[ultra thick,fill=gray!30] (-0.15,-4.15) circle (0.125);
% 	\draw[ultra thick,fill=gray!30] (0.15,-4.15) circle (0.125);
	\draw[fill=gray!30]    (0,-3.5) -- ++(0.3,-0.5) -- ++(-0.6,0) -- ++(0.3,0.5);
	% Triangle
	\draw [ultra thick] (0,-3.5) -- (0.3,-4);
	\draw [ultra thick] (0.5,-4) -- (-0.5,-4);
	\draw [ultra thick] (-0.3,-4) -- (0,-3.5);
	% Dashes
	\draw [ultra thick] (0,-4) -- (-0.2,-4.2);
	\draw [ultra thick] (-0.2,-4) -- (-0.4,-4.2);
	\draw [ultra thick] (0.2,-4) -- (0,-4.2);
	
	\node[draw=none] at (5, -1.5)  (c)     {\large{Minimize weighted sum of}};
	\node[draw=none] at (5, -2)  (c)     {\large{compliance and mass}};
	
	\draw[-latex, line width=1mm,red!80!black] (5,1.25) -- (5,0);
	\node[draw = none] at (5,1.6) () {$2F$}; %{$2F(\xi_p)$};
	
% 	% Rollers
% 	\draw[ultra thick,fill=gray!30] (-0.65,-0.15) circle (0.125);
% 	\draw[ultra thick,fill=gray!30] (-0.65,0.15) circle (0.125);
% 	\draw[fill=gray!30]    (0,0) -- ++(-0.5,0.3) -- ++(0,-0.6) -- ++(0.5,0.3);
% 	% Triangle
% 	\draw [ultra thick] (0,0) -- (-0.5,-0.3);
% 	\draw [ultra thick] (-0.5,-0.3) -- (-0.5,0.3);
% 	\draw [ultra thick] (-0.5,0.3) -- (0,0);
% 	% Rollers
% 	\draw[ultra thick,fill=gray!30] (-0.65,-3.35) circle (0.125);
% 	\draw[ultra thick,fill=gray!30] (-0.65,-3.65) circle (0.125);
% 	\draw[fill=gray!30]    (0,-3.5) -- ++(-0.5,0.3) -- ++(0,-0.6) -- ++(0.5,0.3);
% 	% Triangle
% 	\draw [ultra thick] (0,-3.5) -- (-0.5,-3.8);
% 	\draw [ultra thick] (-0.5,-3.8) -- (-0.5,-3.2);
% 	\draw [ultra thick] (-0.5,-3.2) -- (0,-3.5);
	
% 	% Rollers
% 	\draw[ultra thick,fill=gray!30] (-0.65,-1.90) circle (0.125);
% 	\draw[ultra thick,fill=gray!30] (-0.65,-1.60) circle (0.125);
% 	\draw[fill=gray!30]    (0,-1.75) -- ++(-0.5,0.3) -- ++(0,-0.6) -- ++(0.5,0.3);
% 	% Triangle
% 	\draw [ultra thick] (0,-1.75) -- (-0.5,-2.05);
% 	\draw [ultra thick] (-0.5,-2.05) -- (-0.5,-1.45);
% 	\draw [ultra thick] (-0.5,-1.45) -- (0,-1.75);
	
	\draw[thick,latex-latex] (0,-5) -- (10,-5);
	\node[draw = none] at (5,-4.6) () {$L$};
	\draw[thick] (0,-4.75) -- (0,-5.25);
	\draw[thick] (10,-4.75) -- (10,-5.25);
	
	\draw[thick,latex-latex] (0,0.75) -- (5,0.75);
	\node[draw = none] at (2.5,1.1) () {$L/2$};
	\draw[thick] (0,0.95) -- (0,0.45);
	\draw[thick] (5,0.95) -- (5,0.45);
	
	\draw[thick,latex-latex] (11.5,0) -- (11.5,-3.5);
	\node[draw = none] at (11,-1.75) () {$L/6$};
	\draw[thick] (11.25,0) -- (11.75,0);
	\draw[thick] (11.25,-3.5) -- (11.75,-3.5);
	\end{tikzpicture}
	    \caption{Schematic of the rectangular beam used in Example I.} \label{fig:ex2_schem}
\end{figure}
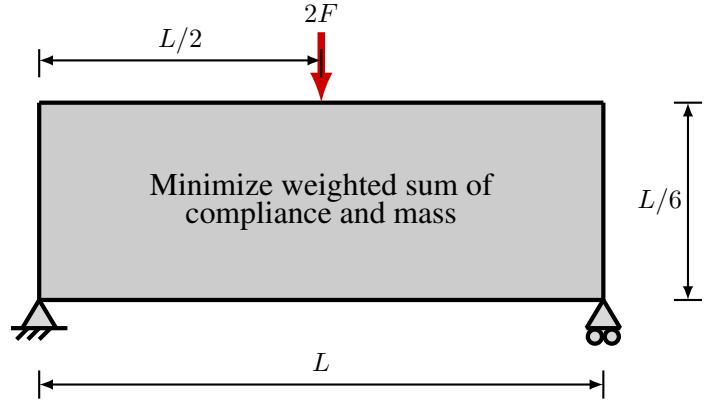

In the first example, the design of a simply supported two-dimensional rectangular beam shown in Figure \ref{fig:ex2_schem} subjected to an uncertain vertical load applied at the midpoint along the top edge is considered. The total length of the beam is assumed as $L = 120$ (in consistent unit). 
Using symmetry condition, only one half of the beam is considered for design. %as reflected in Figure \ref{fig:ex2_schem_right} makes the computations less expensive.
% to minimize an objective function $J(\rhoo)$, {\color{red}which is a weighted combination of the structural compliance $C(\Um;\rhoo,\xii)$ and volume of the structure $V(\rhoo)$} as follows: 
% \begin{equation}
% \begin{split}
%     \min_{\rhoo}~~J(\rhoo) &:=  w_{1}C(\Um;\rhoo,\xii) + w_2 V(\rhoo) + \frac{\kappa_f}{2} \left[\left(\ln P_f -\ln p_a \right)^+\right]^2\\
%         &\text{subject to  }\quad \Km(\rhoo;\xii)\Um(\rhoo;\xii) = F(\xii), 
% \end{split}
% \end{equation} 
% where $\rhoo$ is the density of the elements, $\xii$ denote the uncertain variables, $\Um(\rhoo;\xii)$ is the displacement vector, $w_1$ and $w_2$ are weights for the compliance and volume in the objective, and $\kappa_f$ is the weight for the penalty term for probability of failure $P_f$ exceeding the acceptable failure probability $p_a$. %The compliance is denoted by $C(\rhoo,\Um;\xi)$ which depends on the material density $\rhoo$, the displacement field $\Um$, and the uncertain variable $\xii$. The compliance, volume weight factor used in this example are $w_1$ and $w_2$, respectively. $\kappa_f$ is weight for reliability part. 
The optimization problem is formulated as in \eqref{eq:problem} with $\omega_c = 1$ and $\omega_v = 0.2$. 
% . Refer to Equation \ref{eq:opt_uncon} weights used for compliance and volume are $w_c = 1$ and $w_v = 0.2$, respectively.  %, providing a meaningful comparison between designs that are aware of reliability and those that do not account for uncertainty. 
The vertical load $2F$ is modeled as Gaussian with a mean of one and standard deviation of 0.5.  
% {\color{green}(SD: This distribution is for $F$ or $2F$? Remember, with symmetry, you are only using $F$.)}{\color{green}(MM:  $2F$ for F is used 0.5 for mean and sd of 0.25)}.  
% $2F(\xi_p)=(1+0.5\xi_p)$ where $\xi_p \sim \mathcal{N}(0,1)$, \ie a standard Gaussian random variable with zero mean and unit standard deviation. 
The modulus of elasticity $E_0$ of the material is assumed to follow a lognormal distribution with a mean of one and a standard deviation of $0.1$. The Poisson's ratio $\nu$ follows a uniform distribution between 0.1 and 0.5. These distributions are summarized in Table~\ref{tab:MBB_unc}. 
The failure is defined as the event where the compliance exceeds a specified maximum allowable limit, $C_{\max}=100$. 
Therefore, the reliability constraint is defined as % using the performance metric $g(\rhoo;\xii)$ as 
\begin{equation}
    \Prob (C(\rhoo;\xii) \geq C_{\max}) \leq p_a, 
\end{equation} 
where the target probability of failure is set to $p_a=10^{-2}$. 
% Where the maximum allowed compliance specified is $C_{\max}=100$ to reach the target probability of failure of $p_a=10^{-2}$. 
The penalty parameter for reliability constraint is assumed as $\kappa_f=1500$, which is large enough to implement the probability of failure in this example. 
% Using the symmetry, we only design the right half of the beam as shown in Figure \ref{fig:ex2_schem_right}. 
The design domain is discretized using 4,427 three-noded  triangular elements. 

\begin{table}[!htb]
    \centering
    \begin{tabular}{c|c|c|c} 
    \hline 
       \Tstrut Parameter & Distribution & Mean & Standard deviation \Bstrut \\
        \hline 
        \Tstrut $F$ & Gaussian & 0.5 & 0.25\phantom{0} \Bstrut\\ 
        $E_0$ & Lognormal & 1.0 & 0.1\phantom{00} \Bstrut\\ 
        $\nu$ & Uniform & 0.3 & 0.115 \Bstrut\\
        \hline 
    \end{tabular} 
    \vspace{10pt} 
    \caption{Uncertain parameters and their distributions used in Example I.} 
    \label{tab:MBB_unc}
\end{table} 

\subsubsection{Results} 
% Reliability based topology optimization is conducted in a total number of $12000$ iterations, which is sufficient to ensure convergence for probability of failure and design. At this iteration, we observed that there were no significant changes in the design.
% For comparison, Figure \ref{fig:pth_rob_E1} displays the design obtained by solving the same optimization problem without accounting for reliability (robust design). As it is clear, the reliability-based design (Figure \ref{fig:pth_E1}) is structurally different from the robust design. In RBTO design, the structure has more bars inside to satisfy the probability of failure constraint. 

\begin{figure}[!htb]
\centering
\includegraphics[width=0.7\textwidth]{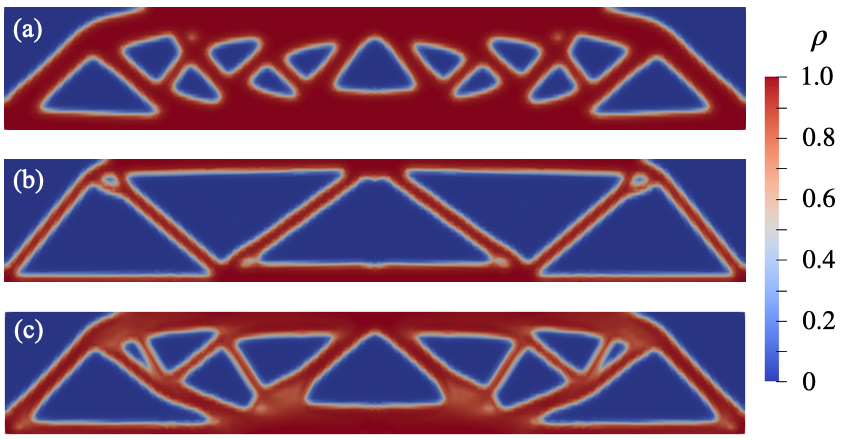}
% \begin{subfigure}[t]{0.84\textwidth} 
%     \centering

%     \begin{subfigure}[t]{\textwidth}
%         \centering
%         \includegraphics[
%             width=\linewidth,
%         ]{Figures/MBB/design_mbb_robust.png}
%         \caption{Robust design of MBB}
%         \label{fig:pth_rob_E1}
%     \end{subfigure}

%     \vspace{0.6em}

%     \begin{subfigure}[t]{\textwidth}
%         \centering
%         \includegraphics[
%             width=\linewidth,
%         ]{Figures/MBB/design_mbb.png}
%         \caption{Reliability-based optimized design}
%         \label{fig:pth_E1}
%     \end{subfigure}
% \end{subfigure}
% \begin{subfigure}[t]{0.14\textwidth}
%     \centering
%     \vspace{-2.6cm}
%     \includegraphics[height=6.0cm]{Figures/pthlegend.png}
%     \label{fig:pthlegend}
% \end{subfigure}

\caption{Comparison of optimized designs obtained from three approaches: (a) proposed LDT based RBTO, (b) robust design without any reliability constraint, and (c) SORA approach with FORM for reliability assessment.}
\label{fig:MBB_design} 
\end{figure}
% We used subset simulation to guide the optimization process toward a reliable design that performs correctly while keeping the computational cost low compared to the Monte Carlo method. 
% However, for the final optimized design, we performed Monte Carlo analysis using $10000$ samples to obtain a more accurate estimate of the probability of failure. The Monte Carlo result yielded $P_f = 0.009$, which closely matches the prescribed probability of failure and falls within the range of the probability of failure obtained using subset simulation. 
% Figure \ref{fig:PF_E1} shows the changes of the estimated probability of failure over the optimization. 

\begin{figure}[!htb]
    \centering
    \includegraphics[width=0.5\linewidth]{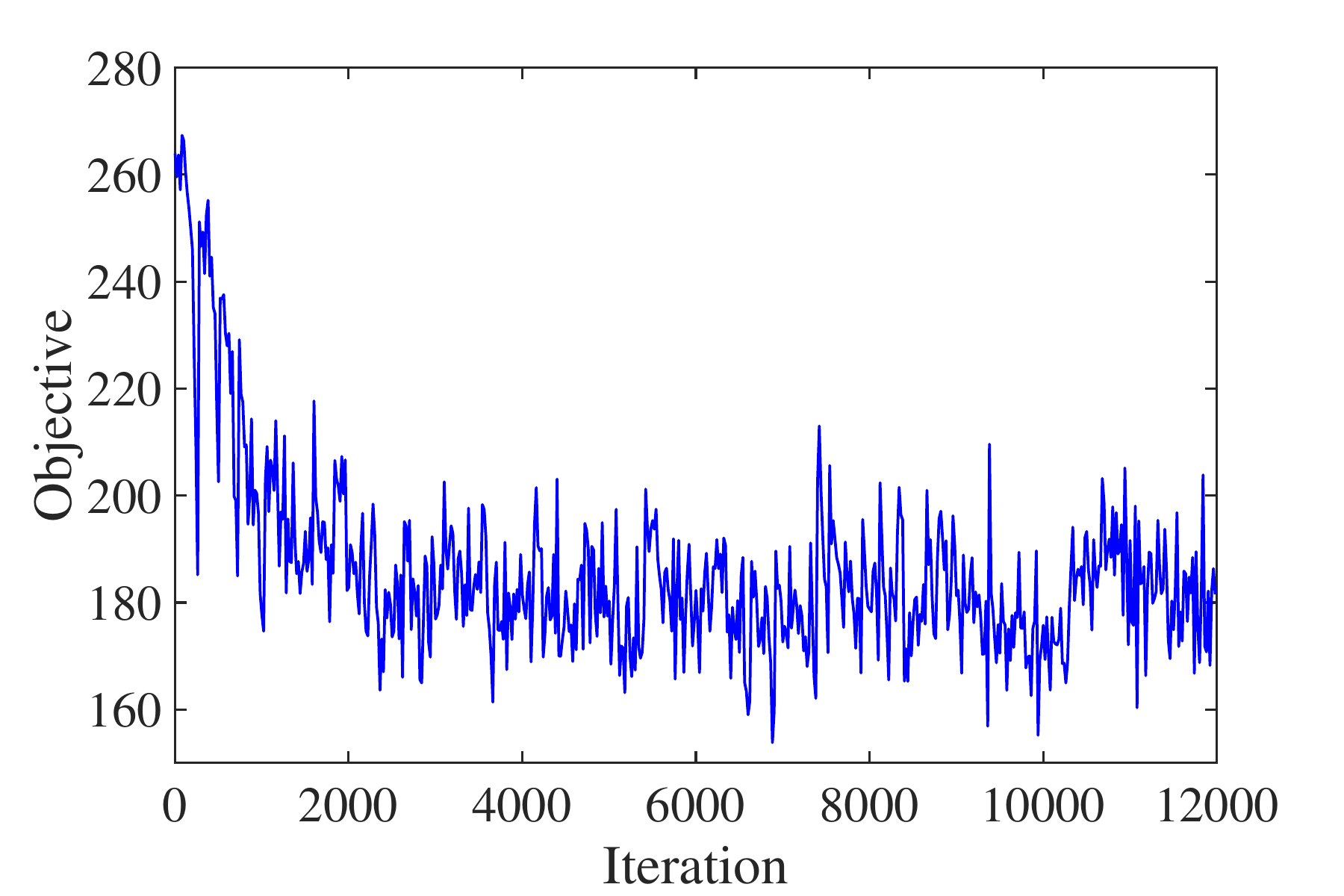}
    \caption{Evolution of objective over the iterations in Example I.}
    \label{fig:Compliance_E1}
\end{figure}

\begin{figure}[!htb]
    \centering
    \includegraphics[width=0.5\linewidth]{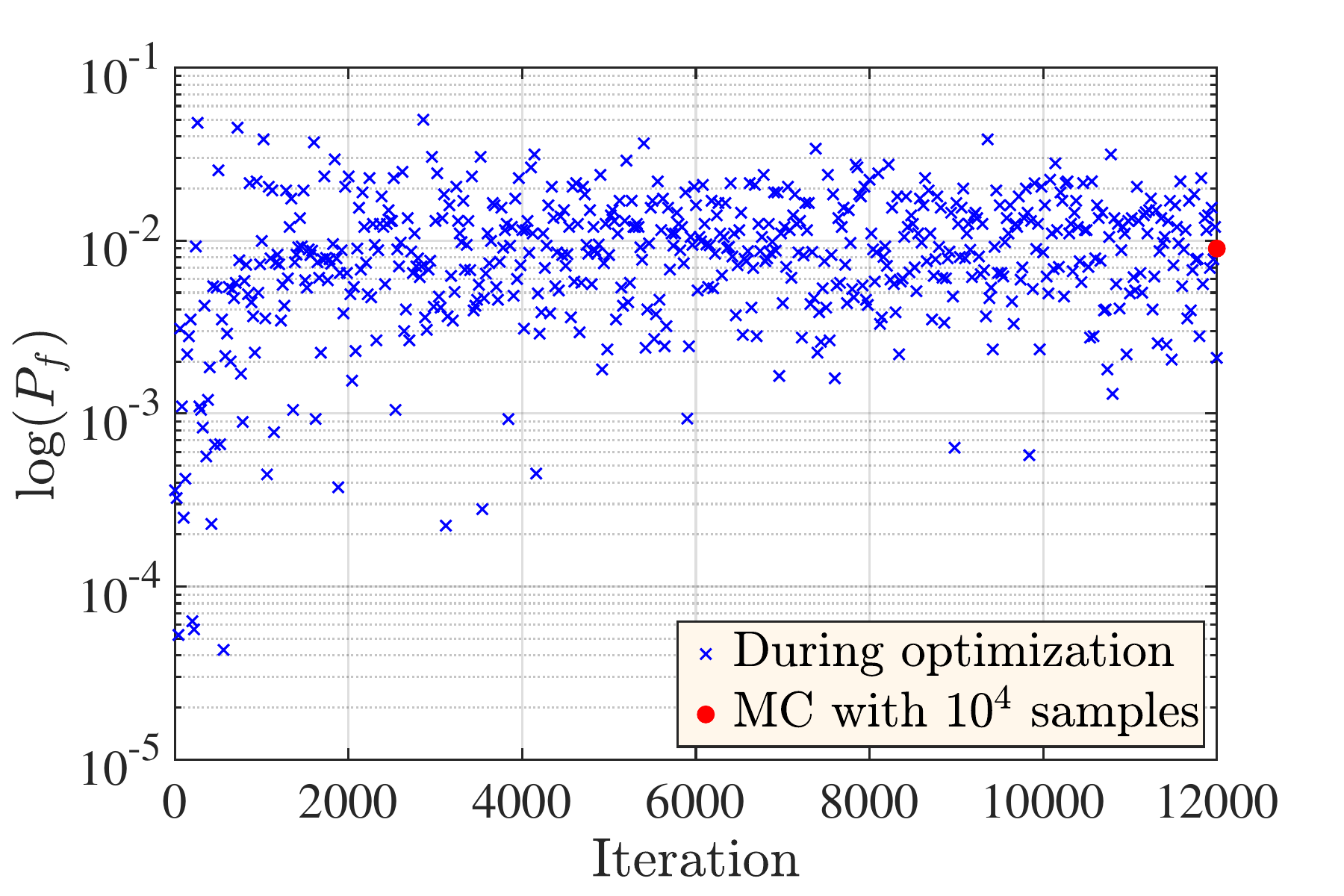}
    \caption{Probability of failure estimated during the optimization along with the probability of failure of the final design estimated using Monte Carlo approach with $10^4$ samples in Example I.}
    \label{fig:PF_E1} 
\end{figure}

Figure~\ref{fig:MBB_design}(a) shows the final design obtained after 12,000 iterations of the proposed optimization framework. The corresponding convergence histories of the objective function and probability of failure are presented in Figures~\ref{fig:Compliance_E1} and~\ref{fig:PF_E1}, respectively. Note that the oscillatory behavior observed in these curves arises from the SGD scheme, which relies on only $n=10$ realizations of the uncertain parameters per iteration. To assess the reliability of the final design from the proposed approach, a Monte Carlo (MC) estimate using 10,000 samples is performed, yielding a probability of failure $P_f = 0.009 < p_a$, thereby satisfying the prescribed failure constraint. 

% The proposed method is implemented over 12,000 iterations. The final design is shown in Figure~\ref{fig:MBB_design}(a). The convergence history of the objective function and probability of failure are shown in Figures~\ref{fig:Compliance_E1} and~\ref{fig:PF_E1}, respectively. Note that the oscillations observed in these plots are due to the use of a SGD based approach which uses only 10 realizations of the uncertain parameters. 
% To validate the final design, a Monte Carlo (MC) based estimation of the probability of failure is performed with 10,000 samples, which shows the final design from the proposed approach has a probability of failure $P_f=0.009<p_a$ satisfying the failure constraint. %While the design obtained from the proposed approach shows a $P_f=0.009<p_a$, the robust design has a $P_f=0.16$, significantly higher than the acceptable value $p_a$. 


When compared to the design obtained from robust optimization (\ie without any reliability constraint), as shown in Figure~\ref{fig:MBB_design}(b), the design from the proposed approach contains more number of bars and each of the bars are thicker. This is due to the fact that the volume used in the design contributes to the objective function, and as a result the optimizer tries to remove as much volume as possible during the optimization. However, by doing so, the probability of failure defined in terms of compliance becomes higher. On the other hand, RBTO explicitly includes the reliability constraints forcing the design to include more volume to satisfy the probability of failure condition. Furthermore, a comparison with the SORA approach \cite{du2004sequential} with first-order reliability method (FORM) for reliability estimation is performed for this problem with the design depicted in Figure~\ref{fig:MBB_design}(c). The failure probability of the final design obtained from this approach is $P_f=0.04$ when estimated using 10,000 MC samples. Note that FORM uses a linear approximation of the failure surface at the MPP, and as a result, the estimation of $P_f$ remains inaccurate for nonlinear failure boundaries leading to a higher $P_f$ of the final design in this case. 
Table~\ref{tab:MBB_result_comp} shows a comparison of estimated $P_f$ using MC and volume used by the final designs obtained from different approaches showing more volume being used to satisfy the reliability constraint. % {\color{red}(SD: Need to add one or two sentences emphasizing the extra volume is needed to satisfy the Pf.)} {\color{blue}(mm:It is concluded that more material is needed to satisfy the probability of failure constraint. The robust design uses less material; however, it is not a reliable design.)}

\begin{table}[!htb]
    \centering
    \begin{tabular}{c|c|c} 
    \hline 
       \Tstrut Method & $P_f$ & Volume fraction \Bstrut \\
        \hline 
        \Tstrut Proposed & 0.009 & 0.65 \Bstrut\\ 
        Robust & 0.16\phantom{0} & 0.31 \Bstrut\\ 
        SORA & 0.04\phantom{0} & 0.48 \Bstrut\\
        \hline 
    \end{tabular} 
    \vspace{10pt} 
    \caption{Comparison of failure probabilities and volume fraction used in the designs obtained from three approaches. 
    % {\color{green}(SD: @MM: Please add the missing values.)}{\color{green}(MM:values are for full mbb and multiplied by weight.)} {\color{green} (SD: Can you check the volume fractions I added in the last column? )}} {\color{green}(SD: I don't understand. Your total volume is $120\times60=7200$. Then $w_v=0.2$ and $w_vV=314.42$. So, $V=1572$ and $V/7200=0.2183$.)}{\color{green}(MM: total volume is $120\times20=2400$. Then $w_v=0.2$ and $w_vV=314.42$. So, $V=1572$ and $V/2400=0.655$.I was mentioned the number of element wrong but the length was correct, I modified that sorry!)
    }
    \label{tab:MBB_result_comp}
\end{table}

\subsection{Example II: Design of a Two-dimensional L-shaped Beam}
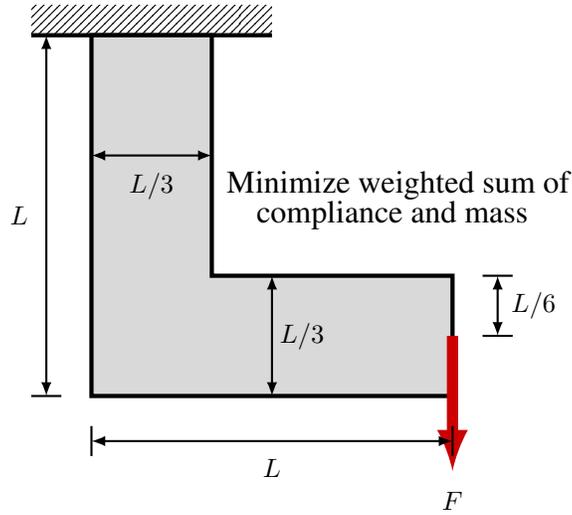
\begin{figure}[!htb]
    \centering
\begin{tikzpicture}[scale=0.8] 
    % \draw[draw=none, fill=gray!50, very thick] (0,2.5) rectangle ++(5,0.25); 
    \draw[ultra thick] (-1,6) -- (3,6);
    \draw[pattern=north east lines,draw=none] (-1,6) rectangle ++(4,0.5); 
    \draw [draw=black,-,ultra thick,fill=gray!30] (0,6) -- (0,0) -- (6,0) -- (6,2) -- (2,2) -- (2,6) -- cycle; 
    \node[draw=none] at (5.1, 3.5)  (c)     {\large{Minimize weighted sum of}};
	\node[draw=none] at (5, 3)  (c)     {\large{compliance and mass}}; 
	\draw[-latex, line width=1.5mm,red!80!black] (6,1) -- (6,-1.25);

	\draw[thick,latex-latex] (-0.75,0) -- (-0.75,6); 
	\node[draw = none] at (-1.2,3) () {$L$}; 
	\draw[thick] (-1,0) -- (-0.5,0); 
% 	\draw[thick] (11.25,-3.5) -- (11.75,-3.5); 

    \draw[thick,latex-latex] (3,0) -- (3,2); 
    \node[draw = none] at (3.5,1) () {$L/3$}; 
    
    \draw[thick,latex-latex] (0,4) -- (2,4); 
    \node[draw = none] at (1,3.5) () {$L/3$}; 
    
    \draw[thick,latex-latex] (0,-0.75) -- (6,-0.75); 
    \node[draw = none] at (3,-1.2) () {$L$}; 
    \draw[thick] (0,-1) -- (0,-0.5); 
    \draw[thick] (6,-1) -- (6,-0.5); 
    
    \draw[thick,latex-latex] (6.75,1) -- (6.75,2); 
    \node[draw = none] at (7.35,1.5) () {$L/6$}; 
    \draw[thick] (6.5,1) -- (7,1); 
    \draw[thick] (6.5,2) -- (7,2); 
    
    \node[draw = none] at (6,-1.75) () {$F$};

    \end{tikzpicture} 
% \includegraphics[scale=1]{figs/Fig7.pdf}
    % \caption{Analytical solution} 
    % \end{subfigure} 
    \caption{A two-dimensional L-shaped beam with uncertain material property and subjected to an uncertain load $F$ at the right edge is used in Example II. } 
    \label{fig:ex3_schem}
\end{figure}

% {\color{green}(SD: Did you use consistent units here? I thought $E$, $\sigma_y$, etc. are in SI units.)}{\color{green}(MM: I used $E$ as GPa, $\sigma_y$ as MPa, etc. we can mention length as mm)} {\color{green}(SD: I have changed the units accordingly. Plz check!)} 
In the second example, we investigate the optimal design of a structure within an L-shaped design domain subjected to material and loading uncertainties with schematic shown in Figure~\ref{fig:ex3_schem}. %depict the two-dimensional configuration of the problem. 
The beam is fixed supported at top and has a length $L=120$ mm. The design domain is discretized into 18,492 three-noded  triangular elements. 
An uncertain vertical load $F$ is applied at the center of the right vertical edge. 
Furthermore, the modulus of elasticity and Poisson's ratio of the material are assumed uncertain. The probability distributions of these uncertain parameters are described in Table~\ref{tab:L_unc}. 
%, where $F$ is assumed as Gaussian distributed with a mean value of 55 and standard deviation of 20. 
% $F(\xi_p)=(55+20\xi_p)$ is applied at the center of the right vertical edge of the domain. Here, $\xi_p \sim \mathcal{N}(0,1)$ is modeled as a standard Gaussian random variable with zero mean and unit standard deviation. 
% The Young’s modulus $E_0$ of the bulk material is assumed to follow a lognormal distribution with mean of $207$ and standard deviation of $20.7$. The Poisson's ratio follows a uniform distribution between $0.1$ and $0.5$. 
% This example is designed to validate the proposed method for the same objective function as example one and Equation \ref{eq:opt_uncon} that combines the expected compliance of the structure, a weighted contribution of the total material volume, and a penalty term proportional to the probability of failure.
\begin{table}[!htb]
    \centering
    \begin{tabular}{c|c|c|c|c} 
    \hline 
       \Tstrut Parameter & Unit & Distribution & Mean & Std dev \Bstrut \\
        \hline 
        \Tstrut $F$ & N & Gaussian & \phantom{0}55.0 & \phantom{0}20.0\phantom{00} \Bstrut\\ 
        $E_0$ & kN/mm$^2$ & Lognormal & 207.0 & \phantom{0}20.7\phantom{00} \Bstrut\\ 
        $\nu$ & --- & Uniform & \phantom{00}0.3 & \phantom{00}0.115 \Bstrut\\
        \hline 
    \end{tabular} 
    \vspace{10pt} 
    \caption{Uncertain parameters and their distributions used in Example II.} 
    \label{tab:L_unc}
\end{table}

The optimization problem is defined with a compliance weight of $\omega_c=1$ and a volume weight of $\omega_v=10$ in \eqref{eq:problem}. 
% , and a reliability penalty weight of $\kappa_f=10^4$. 
As the L-shaped domain exhibits stress concentration at the inner corner, it is essential to track the stress field during the optimization process. To emphasize this critical aspect, the performance function in the reliability constraint is modified from the compliance based formulation used in Example I to a $p$-norm stress-based one, which provides more control over stresses developed in the structure. 
Therefore, the failure is defined as the event where the $p$-norm stress from \eqref{eq:pnorm}, an approximation of the maximum stress, exceeds an allowable limit of
$\sigma_\mathrm{yield}=370$ MPa, 
which corresponds to the yield stress of AISI 1018 cold-drawn steel. Hence, in this example, the reliability constraint is 
\begin{equation}
    \Prob (\sigmap\geq \sigma_\mathrm{yield})\leq p_a, 
\end{equation}
% material properties of the structural steel. 
where the allowable probability of failure is set to $p_a=10^{-3}$, and the penalty parameter for reliability constraint is assumed as $\kappa_f=10^4$. Since a penalty-based approach is used to implement the violation of the reliability constraint, and no penalty is applied when estimated $P_f<p_a$, large oscillations can be observed in $P_f$ with a gradient-based approach. Therefore, in the remaining examples, a conservative threshold of $p_a/2$ is employed in Algorithm~\ref{alg:ldt-rbto} to yield a safer design and dampen oscillations near $p_a$.

\subsubsection{Results}
\begin{figure}[!htb]
\centering
\includegraphics[width=0.7\textwidth]{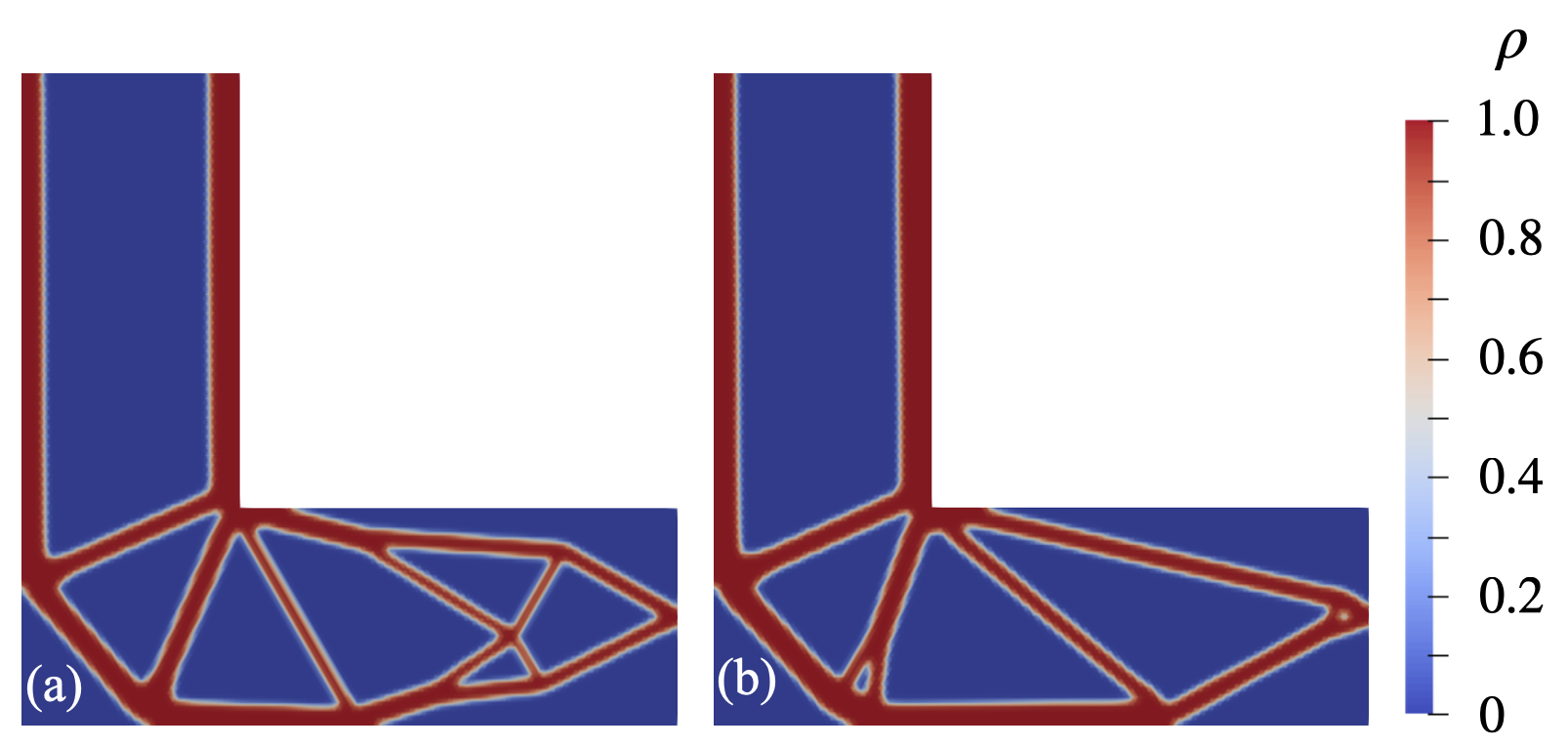}
% % ---------- (a) ----------
% \begin{subfigure}[t]{0.40\textwidth}
%     \centering
%     \includegraphics[width=\linewidth]{Figures/E2/robust optimized lshape beam.png}
%     \caption{Robust optimized topology of L-shape beam}
%     \label{fig:robust_E2}
% \end{subfigure}
% \hfill
% % ---------- (b) ----------
% \begin{subfigure}[t]{0.40\textwidth}
%     \centering
%     \includegraphics[width=\linewidth]{Figures/E2/Optimised L-Shape beam.png}
%     \caption{Reliability-based optimized design}
%     \label{fig:reliability_E2}
% \end{subfigure}
% \hfill
% % ---------- Legend ----------
% \begin{subfigure}[t]{0.14\textwidth}
%     \centering
%     \includegraphics[height=6.8cm]{Figures/pthlegend.png}
%     \label{fig:pthlegend}
% \end{subfigure}

\caption{Comparison of optimized designs obtained from two approaches: (a) proposed LDT based RBTO and (b) robust design without any reliability constraint.}
\label{fig:L-Shape_design}
\end{figure}

Figure \ref{fig:L-Shape_design}(a) demonstrates the RBTO optimized design obtained using the stress driven failure criterion with a volume fraction of 0.28, while
Figure~\ref{fig:L-Shape_design}(b) presents the robust design (\ie without the reliability constraint) with a volume fraction of 0.26. 
A comparison between these two designs 
% RBTO and robust design 
reveals that the optimizer added more bars inside the design domain to satisfy the stress-based failure condition in the reliability constraint. 
In particular, more inclined bars are added near the L-shaped beam corner and the load location to distribute the developed stress within the design domain, whereas the vertical part of the design remains similar. %{\color{green}(SD: Please add the volume fractions for each of these designs.)}
% Also, the structure is divided to smaller voids near the applied force location. The vertical part of the beam remains similar in both designs.

\begin{figure}[!htb]
    \centering
    \includegraphics[width=0.5\linewidth]{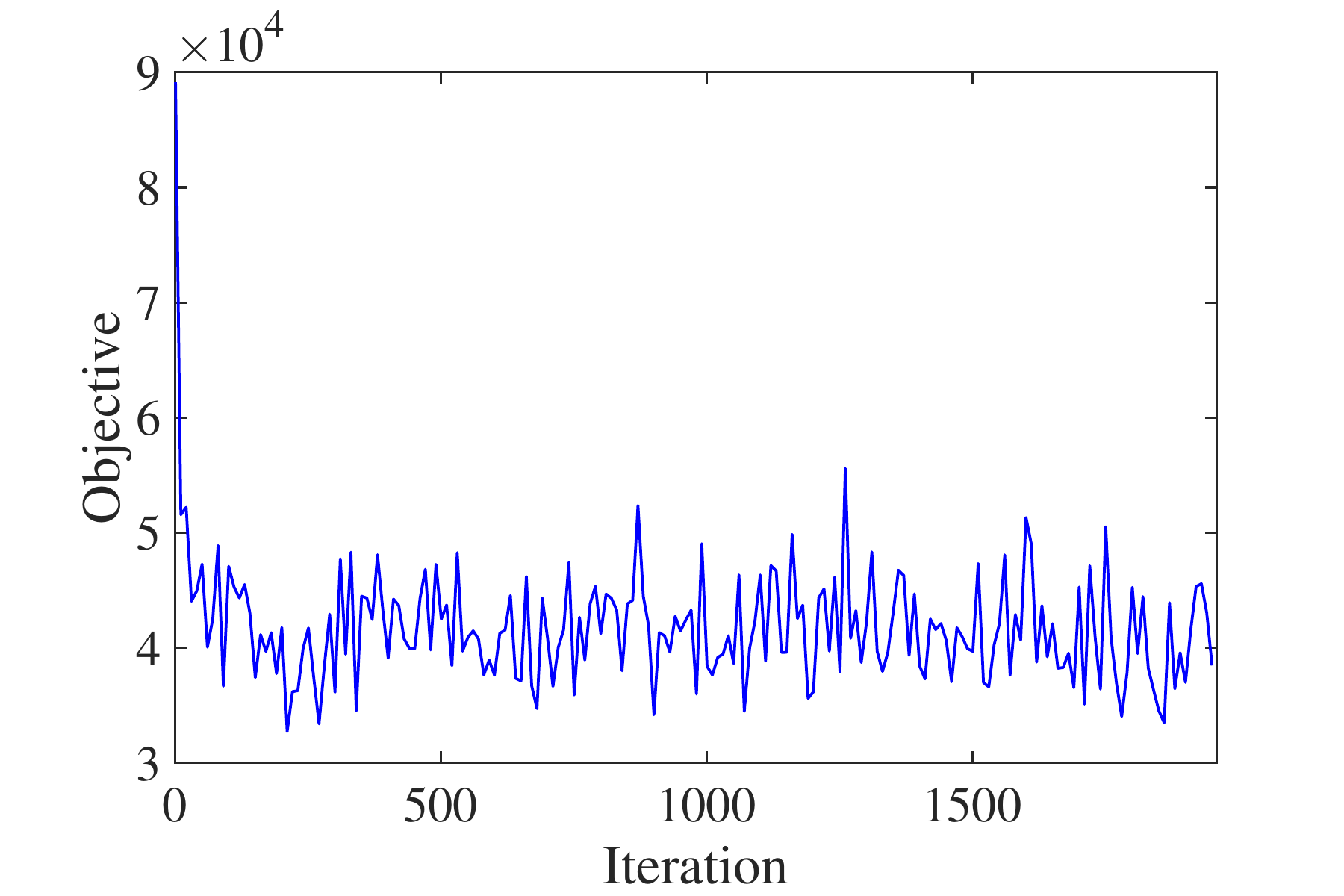}
    \caption{Evolution of objective over the iterations in Example II.}
    \label{fig:Compliance_E2}
\end{figure}

\begin{figure}[!htb]
    \centering
    \includegraphics[width=0.5\linewidth]{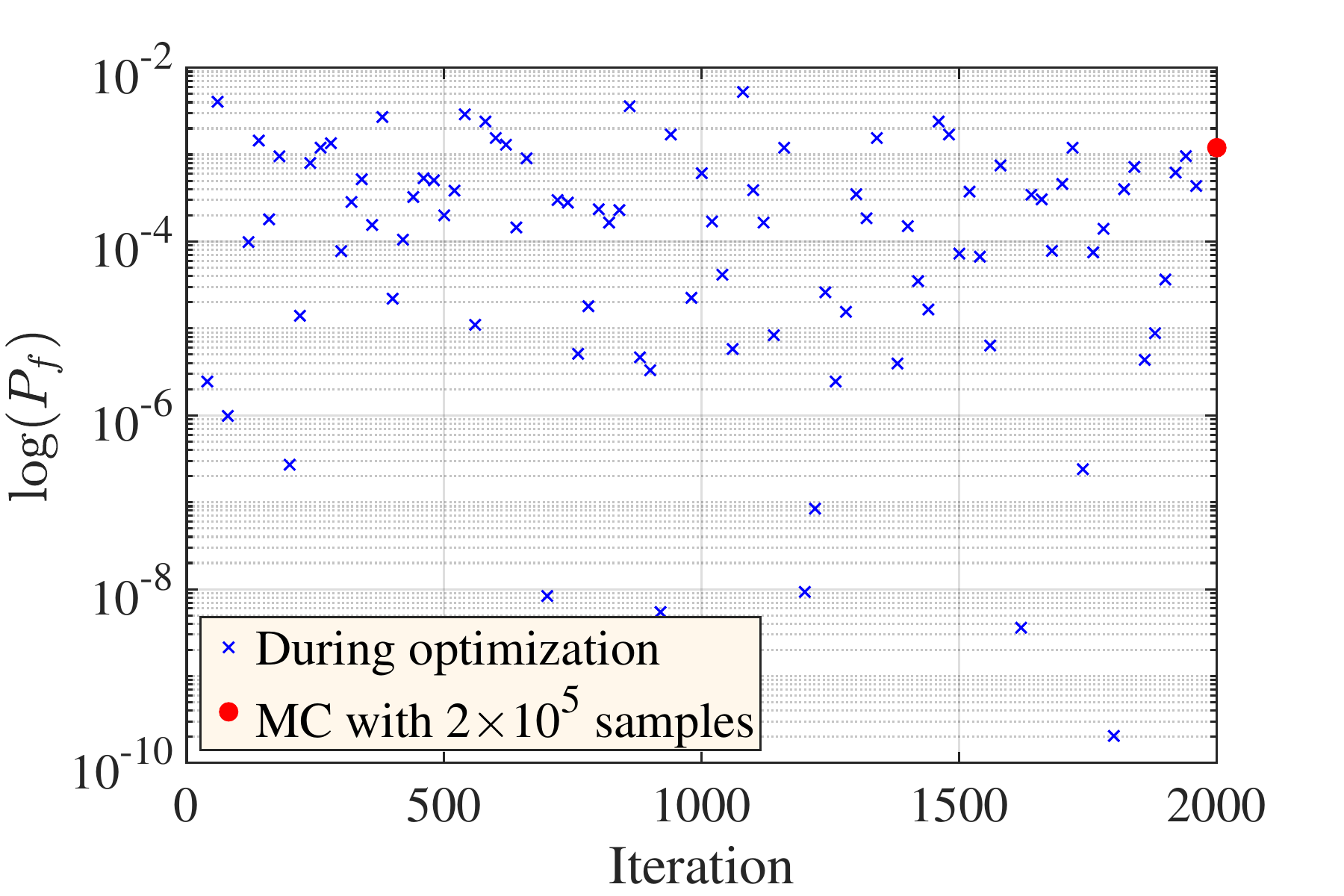}
    \caption{Probability of failure during optimization along with the probability of failure of the final design estimated using Monte Carlo approach in Example II.}
    \label{fig:pf_E2}
\end{figure}

\begin{figure}[!htb]
\centering
\includegraphics[width=0.5\linewidth]{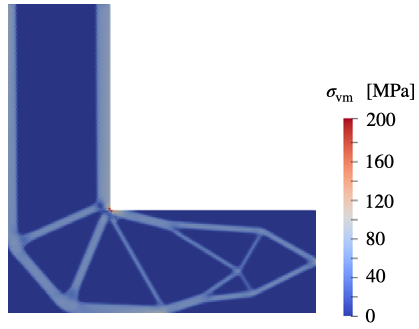}
% \setlength{\tabcolsep}{0pt}

% \begin{tabular}{@{}c@{}c@{}}
% \includegraphics[width=0.48\textwidth]{Figures/E2/Stress Distribution L-Shape.png}%
% &
% \includegraphics[height=8.8cm]{Figures/E2/Stress L-Shape legend.png}

% \end{tabular}

\caption{Typical stress distribution inside the L-shaped beam in Example II.}
\label{fig:Stress_E2}
\end{figure}

% Figure \ref{fig:pf_E2} depicts the convergence of the estimated probability of failure during optimization with the 
Figures~\ref{fig:Compliance_E2} and \ref{fig:pf_E2} depict the convergence of the objective and the probability of failure during optimization, respectively, illustrating that after some initial fluctuations, these quantities stabilizes with $P_f$ near the allowable limit $p_a$.  
% probability of failure stabilizes at the prescribed target value. 
% It should be noted that, to obtain a safer design in the remaining examples, $p_a/2$ was adopted as the constraint violation threshold.  
When evaluated via MC simulation using $2\times10^5$ samples, the probability of failure of the proposed RBTO design yields $P_f = 0.0012$. In comparison, the robust design yields $P_f = 0.0158$, approximately 13th times higher than that of the RBTO design showing the usefulness of incorporating the reliability constraint in the optimization problem while using similar volume fraction. Figure~\ref{fig:Stress_E2} illustrates the stress distribution within the final optimized structure for a representative case. Hence, this example shows the efficacy of the proposed method in accounting for critical stresses developed in the structure by incorporating a $p$-norm-based maximum stress measure into the reliability constraint. 

% incorporation of the $p$-norm stress within the reliability constraint effectively consider critical stress while maintaining overall structural reliability. 
% On the other hand, 
% the corresponding probability of failure of the RBTO structure, assessed 
% under identical conditions, is $P_f = 0.0012$. As is evident, the
% probability of failure of the robust design is approximately $13$ times
% higher than that of the RBTO design.

% oscillation of the  probability of failure, illustrating that, after some fluctuations, the probability of failure remains in the prescribed value of probability of failure. It should be mentioned that to get a safer design in remaining examples we used $p_a/2$ as the constraint violation. As we are comparing the robust and reliable design, the probability of failure of the robust design were determined by Monte Carlo simulation using $200000$ number of samples as $P_f = 0.0158$. We should emphasis on the probability of failure of the RBTO structure using the same number of samples with Monte Carlo which is $P_f = 0.0012$. As it is clear, the robust design probability of failure is approximately $13$ times higher than that of the RBTO design. 

% The convergence plot of compliance and volume summation in Figure \ref{fig:Compliance_E2} confirms the reliable design. The summation over volume and compliance is decreasing which is the goal of topology optimization in this case as the minimization of volume and compliance.

% Furthermore, Figure \ref{fig:Stress_E2} illustrates the stress distribution within the final optimized structure. The results clearly show that the design successfully distributes material to reduce stress concentration near the corner region particularly prone to high stress under the applied load. This demonstrates that incorporating the $p$-norm stress within the reliability condition effectively consider critical stress while maintaining overall structural reliability. 

\subsection{Example III: Design of a Three-dimensional Cantilever Beam} 
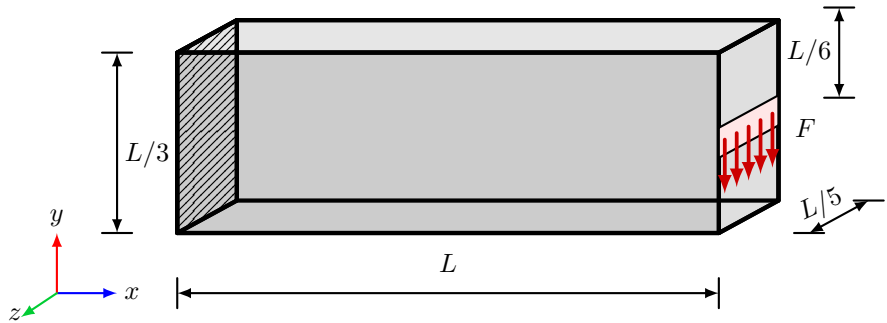
\begin{figure}
    \centering
% 3D Cantilever Beam with Distributed Load
% Geometry: Length = L, Height = L/3, Width = L/5

\begin{tikzpicture}[
  scale=0.8,
  line join=round, line cap=round,
  x={(1cm,0cm)}, y={(0cm,1cm)}, z={(0.55cm,0.30cm)}
]

% -----------------------
% Geometric parameter
% -----------------------
\def\L{9}   % Choose numerical value for visualization

\pgfmathsetmacro{\H}{\L/3}
\pgfmathsetmacro{\W}{\L/5}

% Distributed load patch size
\pgfmathsetmacro{\patchH}{\H/6}
\pgfmathsetmacro{\patchW}{\W*0.8}
\def\nArrows{5}

% -----------------------
% Coordinates
% -----------------------
\coordinate (A) at (0,0,0);
\coordinate (B) at (\L,0,0);
\coordinate (C) at (\L,-\H,0);
\coordinate (D) at (0,-\H,0);

\coordinate (A2) at (0,0,\W);
\coordinate (B2) at (\L,0,\W);
\coordinate (C2) at (\L,-\H,\W);
\coordinate (D2) at (0,-\H,\W);

% -----------------------
% Beam (3D cuboid)
% -----------------------
\draw[ultra thick, fill=gray!15] (A2) -- (B2) -- (C2) -- (D2) -- cycle;
\draw[ultra thick, fill=gray!25] (B) -- (B2) -- (C2) -- (C) -- cycle;
\draw[ultra thick, fill=gray!20] (A) -- (B) -- (B2) -- (A2) -- cycle;
\draw[ultra thick, fill=gray!40] (A) -- (B) -- (C) -- (D) -- cycle;

\draw[ultra thick] (A) -- (B) -- (C) -- (D) -- cycle;
\draw[ultra thick] (A) -- (A2) (B) -- (B2) (C) -- (C2) (D) -- (D2);
\draw[ultra thick] (A2) -- (B2) -- (C2) -- (D2) -- cycle;

% -----------------------
% Fixed boundary at x=0
% -----------------------
\draw[ultra thick] (A) -- (D) -- (D2) -- (A2) -- cycle;

% \draw[thick, fill=gray!10]
%   (0,0,0) -- (0,-\H,0) -- (0,-\H,\W) -- (0,0,\W) -- cycle;

% --- Full rectangular support wall with hatch pattern
\draw[
    thick,
    fill=gray!10,
    pattern=north east lines,
    pattern color=black
]
(0,0,0) --
(0,-\H,0) --
(0,-\H,\W) --
(0,0,\W) --
cycle;

% \node[anchor=east] at (-0.8,-\H/2,\W/2) {\small Fixed};

% -----------------------
% Distributed load at mid-height of right face
% -----------------------
\pgfmathsetmacro{\yTop}{-\H/2 + \patchH/2}
\pgfmathsetmacro{\yBot}{-\H/2 - \patchH/2}
\pgfmathsetmacro{\zLeft}{0}
\pgfmathsetmacro{\zRight}{\W}

% Load patch
\draw[thick, fill=red!10]
  (\L,\yTop,\zLeft) -- (\L,\yTop,\zRight) --
  (\L,\yBot,\zRight) -- (\L,\yBot,\zLeft) -- cycle;

\pgfmathsetmacro{\yMid}{(\yTop+\yBot)/2}

\foreach \i in {0,...,\numexpr\nArrows-1\relax}{
  \pgfmathsetmacro{\zz}{\zLeft + (\i+0.5)*((\zRight-\zLeft)/\nArrows)}
  \draw[-latex, line width=0.5mm, red!80!black]
    (\L,\yMid,\zz) -- (\L,\yMid-0.9,\zz);
}

\node at (\L+1,-\H/2,\W/2-0.1) {$F$};

% -----------------------
% Dimensions
% -----------------------
% Length
\draw[thick, latex-latex] (0,-\H-1,0) -- (\L,-\H-1,0); 
\draw[thick] (0,-\H-0.75,0) -- (0,-\H-1.25,0); 
\draw[thick] (\L,-\H-0.75,0) -- (\L,-\H-1.25,0); 
\node at (\L/2,-\H-0.5,0) {$L$}; 

% Height
\draw[thick, latex-latex] (-1,0,0) -- (-1,-\H,0); 
\draw[thick] (-1.25,0,0) -- (-0.75,0,0); 
\draw[thick] (-1.25,-\H,0) -- (-0.75,-\H,0); 
\node at (-0.5,-\H/2-0.2,0) {$L/3$}; 

% Height
\draw[thick, latex-latex] (11,0.75,0) -- (11,-\H/2+0.75,0); 
\draw[thick] (10.75,0.75,0) -- (11.25,0.75,0); 
\draw[thick] (10.75,-\H/2+0.75,0) -- (11.25,-\H/2+0.75,0); 
\node at (10.5,-\H/2+1.5,0) {$L/6$}; 

%width
\draw[thick, latex-latex]
    (\L+1.5,-\H,0) -- (\L+1.5,-\H,\W); 

\node[rotate=25] at (\L+1.2,-\H+0.2,\W/2) {$L/5$};
\draw[thick] (\L+1.25,-\H,0) -- (\L+1.75,-\H,0); 
\draw[thick] (\L+1.25,-\H,\W) -- (\L+1.75,-\H,\W); 

% coordinate axes
\draw[thick,-latex,blue] (-2,-4,0) -- (-1,-4,0); 
\draw[thick,-latex,red] (-2,-4,0) -- (-2,-3,0); 
\draw[thick,-latex,green!80!blue] (-2,-4,0) -- (-2.6,-4.4,0); 

\node[] at (-0.75,-4,0) {$x$}; 
\node[] at (-2,-2.75,0) {$y$}; 
\node[] at (-2.7,-4.35,0) {$z$}; 

\end{tikzpicture}
\caption{3D cantilever beam used in Example III.}
    \label{fig:3d_Cantilever_beam}
\end{figure}

In the final example, we use a three-dimensional cantilever beam supported at the left end and an uncertain vertical line load $F$ applied at the middle of the right end as shown in Figure~\ref{fig:3d_Cantilever_beam}. With consistent units, the length of the beam is set to $L = 60$. The design domain is discretized using 66,490 tetrahedral finite elements to capture the three-dimensional response.
% Finally, we applied the proposed reliability-based topology optimization (RBTO) method to the design of a three-dimensional cantilever beam. The model has been shown in Figure \ref{fig:Cantilever beam}. The beam is fully supported at one end. The length of the beam was set to $L = 60$. The design domain is discretized using $66490$ tetrahedral finite elements to capture the three-dimensional response.\\ 
% The structure is subjected to an uncertain vertical load of 
% The magnitude of the uncertain load $F$ applied at the midpoint along the right edge. 
The load $F$ is modeled as Gaussian distributed with a mean of 0.5 and standard deviation of 0.5. The distributions of the modulus of elasticity and Poisson's ratio are assumed to be the same as used in Example I. 
% $F=(0.5+0.5\xi_p)$, where $\xi_p \sim \mathcal{N}(0,1)$ is a standard Gaussian random variable. 
% The material’s Young’s modulus is assumed to follow a lognormal distribution with a mean of $1$ and a standard deviation of $0.5$. The Poisson's ratio considered as a uniform distribution between $0.1$ to $0.5$. 
In this example, the failure criterion is defined as the compliance exceeding $C_{\max}=5\times10^4$, with the maximum allowable probability of failure specified as $p_a=10^{-2}$. The weights used for compliance and volume in \eqref{eq:problem} are $\omega_c = 1$ and $\omega_v = 10$, respectively. The penalty parameter used to enforce the reliability constraint is assumed as $\kappa_f = 10^4$. 

% 3D-style Cantilever Beam with point load at mid-height of RIGHT edge (front face)

\subsubsection{Results}
% {\color{green}The design parameter $\rho$ are updated over a total of $560$ iterations, the optimization stopped there since there were no considerable changes in design, the probability of failure oscillations are less than $4\%$ in the last $60$ iterations and the target probability of failure were achieved. (SD: Did you use the same criterion for all the examples? If not, what are they for Ex I and II?) }{\color{green}(I added the stopping criteria in first paragraph of numerical solution. we can remove this part)}

\begin{figure}[!htb]
\centering
\includegraphics[width=0.45\linewidth]{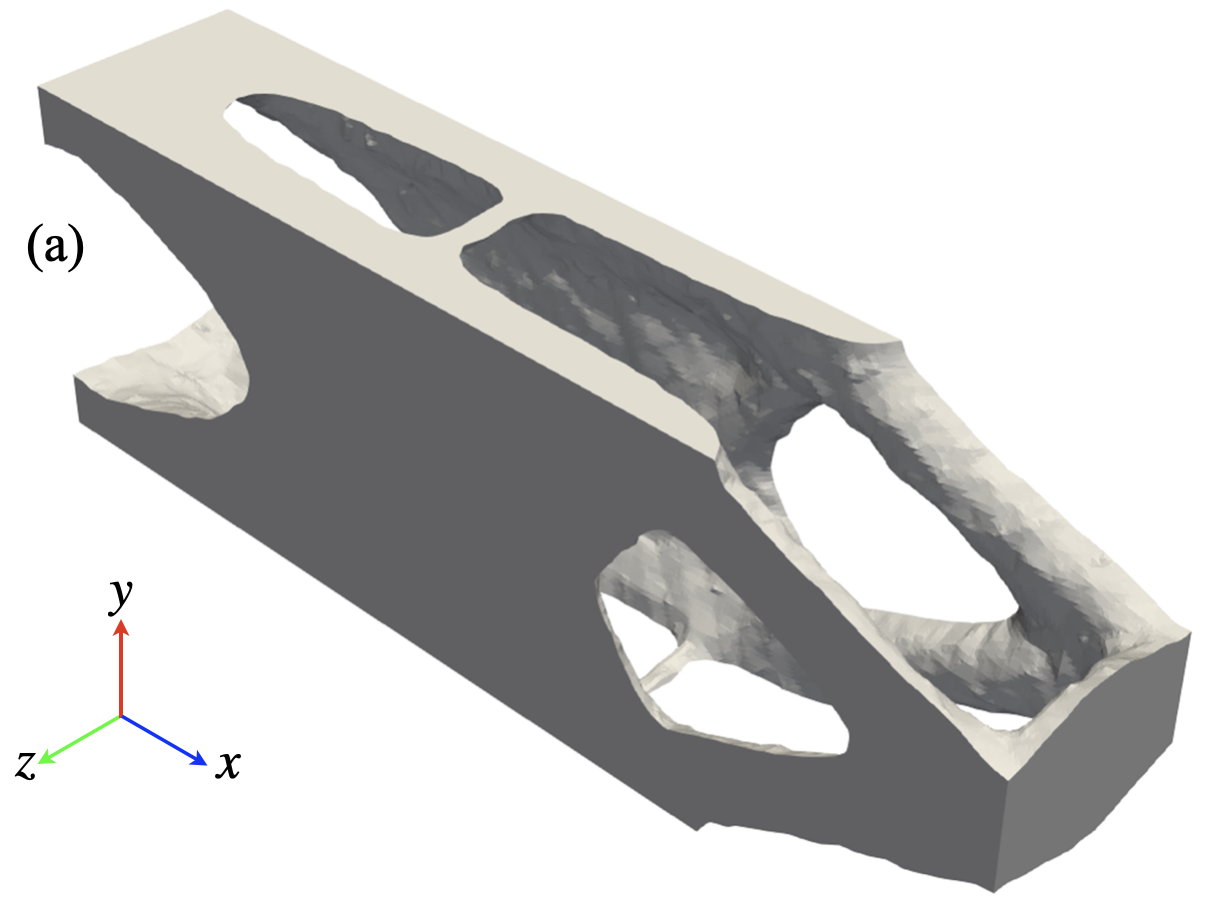}~~\includegraphics[width=0.45\linewidth]{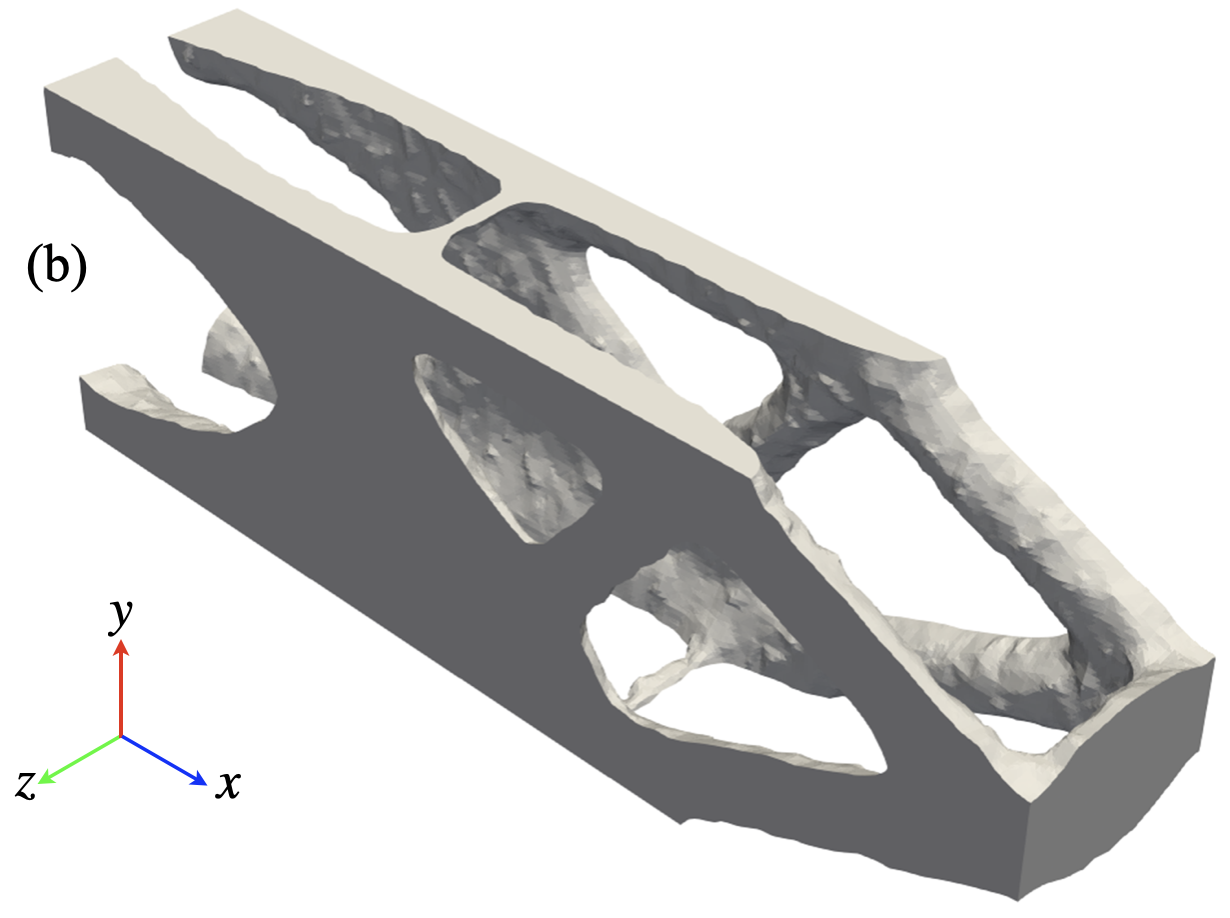}
% % ---------- (a) ----------
% \begin{subfigure}[t]{0.45\textwidth}
%     \centering
%     \includegraphics[width=\linewidth]{Figures/3d/3D_RBTO.png}
%     \caption{Robust optimized topology of 3D cantilever beam}
%     \label{fig:robust_E3}
% \end{subfigure}
% \hfill
% % ---------- (b) ----------
% \begin{subfigure}[t]{0.45\textwidth}
%     \centering
%     \includegraphics[width=\linewidth]{Figures/3d/3D_robust.png}
%     \caption{Reliability-based optimized 3D cantilever beam}
%     \label{fig:reliability_E3}
% \end{subfigure}
% \hfill

\caption{Comparison of optimized 3D cantilever beam using (a) the proposed RBTO approach and (b) robust design without any reliability constraint.}
\label{fig:3d_design}
\end{figure}

\begin{figure}[!htb]
    \centering
    \includegraphics[width=0.5\linewidth]{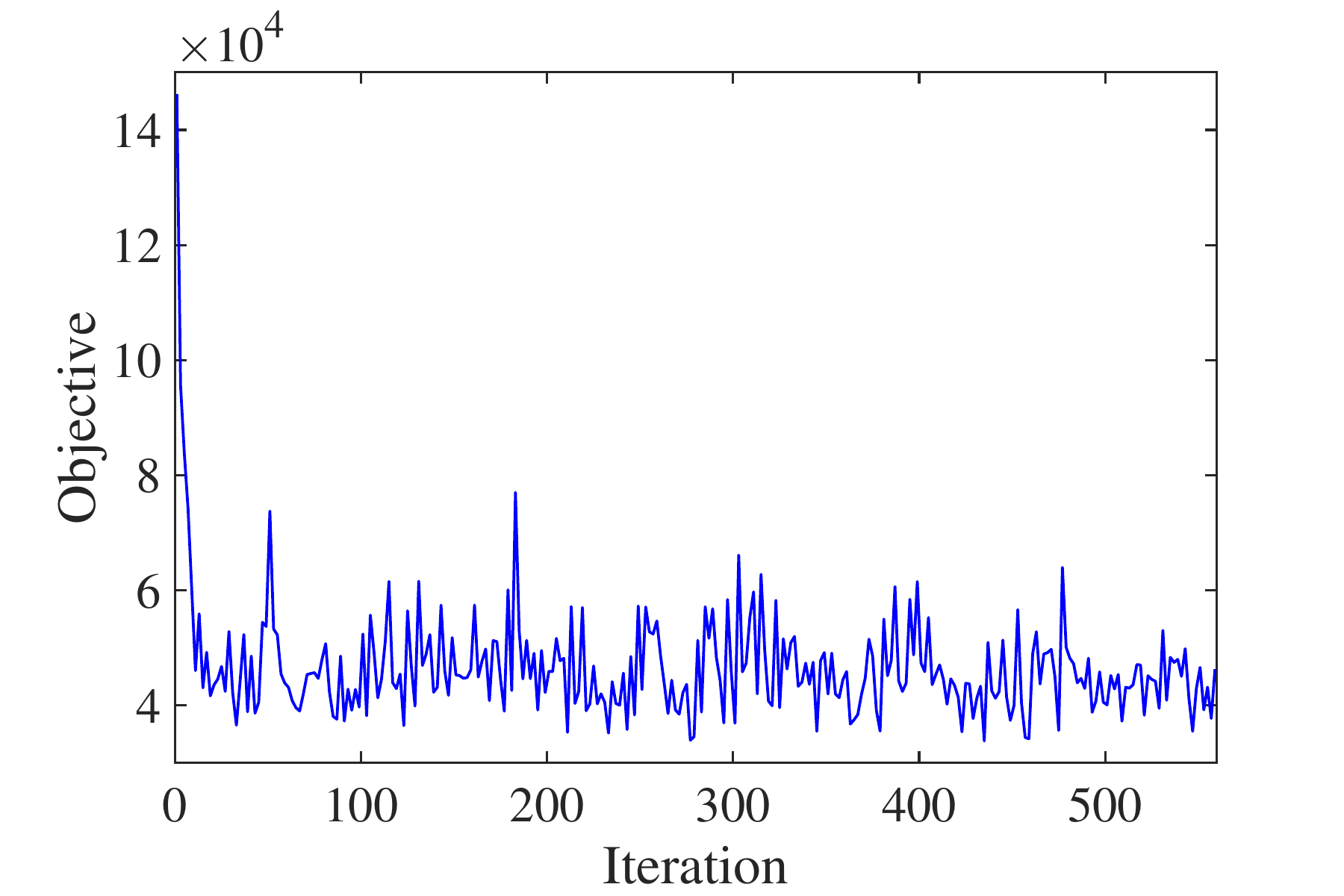}
    \caption{Evolution of objective over the iterations in Example III.}
    \label{fig:Compliance_E3}
\end{figure}

\begin{figure}[!htb]
    \centering
    \includegraphics[width=0.5\linewidth]{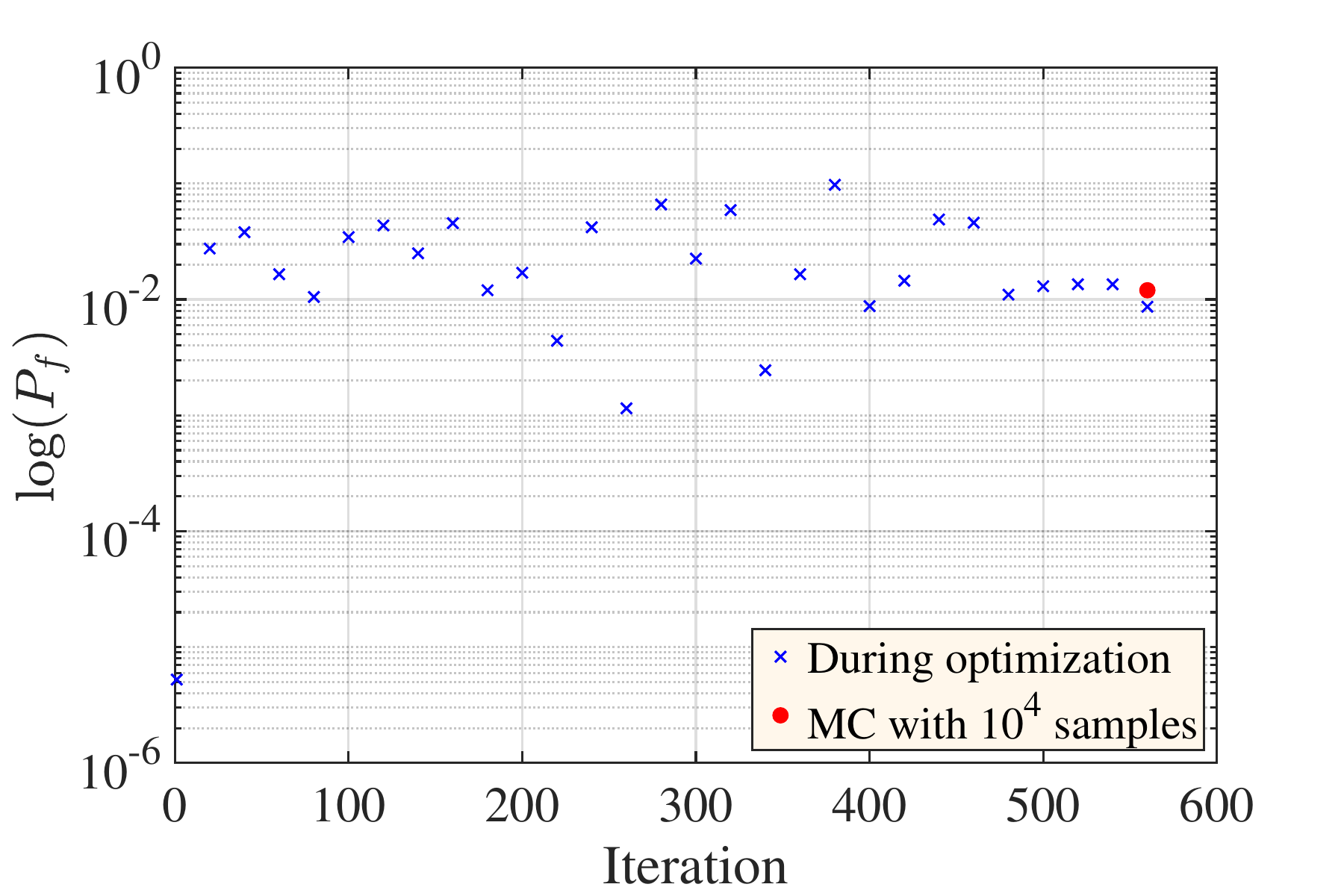}
    \caption{Probability of failure during optimization  along with the probability of failure of the final design estimated using Monte Carlo approach in Example III.}
    \label{fig:pf3d}
\end{figure}

Figure~\ref{fig:3d_design}(a) depicts the final design obtained from the proposed RBTO approach with a volume fraction of 0.23. For comparison, Figure~\ref{fig:3d_design}(b) shows the robust design obtained without considering the reliability constraint, which uses a volume fraction of 0.19. The RBTO design here shows thicker elements and better formed webs on two sides compared to the robust design. The use of more material in the RBTO design is justified by the fact that more material is needed to satisfy the reliability constraint. 
The evolution of objective during the optimization is presented in Figure~\ref{fig:Compliance_E3}, which shows smaller oscillations after 500 iterations. Similarly, Figure~\ref{fig:pf3d} depicts the probability of failure during the optimization. It starts from a small value, as the design is initialized with density of all elements to one, and showing convergence after 500 iterations terminating the optimization. When the final designs are evaluated using an MC approach with 10,000 samples, the RBTO design has a probability of failure of 0.012, slightly higher than $p_a=0.01$, whereas the robust design has a probability of failure of 0.047, four times the probability of failure of the proposed RBTO design. 
Since the evaluation of the performance function $g(\rhoo;\xii)$ and the failure probability are more computationally intensive in three-dimensional problems, the proposed stochastic gradient-based RBTO method offers a significant computational savings compared to other sampling-based methods while achieving a final design satisfying the reliability constraint. 

% by reducing the number of evaluations required for accurate probability of failure estimation in other approaches, such as the Monte Carlo method. 

% Figure \ref{fig:3d_design}(a) shows the RBTO design obtained from the proposed approach. For comparison, Figure \ref{fig:robust_E3} presents the corresponding  robust design. A comparison of these two designs reveals that the RBTO solution have complete bars at upper and lower left edge at fixed support which is not appeared in robust design. These full bar allows the optimizer to satisfy the reliability constraint. It is shown that the robust design have more voids in comparison to the RBTO design. It is remarkable that in RBTO design, the web and flange structure has appeared. 
% The illustrated designs show cells with more than $40\%$ material volume fraction and are smoothed using a mesh refined to four times the original number of elements to have a smoother structure. 

% Figure \ref{fig:pf3d} depicts the oscillations of the probability of failure. It exhibits initial fluctuations but converge after $560$ iterations. The value of probability of failure starts from a small value as the domain is fully filled with material at the beginning and it is increased while the optimizer removes material. The probability of failure remains approximately constant in last $100$ iterations indicating convergence. Both RBTO and robust design were analyzed using Monte Carlo by $10000$ samples to evaluate their probability of failure, which shows $0.012$ and $0.047$ respectively. These results confirm the reliability of RBTO design with respect to robust design. 
% The changes in summation of compliance and volume values is presented in Figure \ref{fig:Compliance_E3} which shows not drastic changes in this quantity near terminated iteration. 
% It is important to note that the evaluation of the limit state function and the failure probability are more computationally intensive in three-dimensional problems. The proposed stochastic gradient–based RBTO method, combined with subset simulation, offers a computationally efficient alternative by reducing the number of evaluations required for accurate probability of failure estimation in other approaches, such as the Monte Carlo method. 

\section{Conclusions} \label{sec:conclusions}
A reliability-based topology optimization framework is introduced herein that combines large deviation theory with a stochastic gradient descent-based optimization approach to alleviate the tremendous computational cost involved with RBTO. In particular, the proposed method uses the large deviation theory to estimate the gradient of the probability of failure, which enables the use of an efficient mini batch  
% ({\color{green} MMR: meaning of "efficient mini batch" should be made clear in the numerical sections before writing it in the conclusion?})
stochastic gradient descent method. This allows effective optimization while incorporating uncertainty using a small number of samples per iteration making the approach practical for large scale RBTO problems. 
% The main contribution of this work is an efficient way to estimate the probability of failure and its gradient using only a small number of samples. 
% This is especially practical for RBTO problems where traditional sampling methods become computationally expensive. By using large deviation theory, the proposed approach addressed the related issues in other methods such as Monte Carlo. 
Three numerical examples are used to demonstrate the efficacy of the proposed approach. Uncertainties are incorporated in material properties and loading conditions using known probability distributions. The first example illustrates the RBTO design of a 2D rectangular beam that results in higher number of thicker bars inside the domain compared to the robust design while satisfying the accepted probability of failure of $0.01$. The second example incorporates a stress-based failure criterion in an L-shaped beam. The proposed method effectively considers the stress distribution and adds more bars compare to robust design to satisfy the accepted probability of failure of $0.001$. The third example implements the approach for the design of a 3D cantilever beam 
% which requires more analyses time and elements. However, the RBTO design with target probability of failure of $0.01$ achieved successfully. 
resulting in a RBTO design that shows webs forming at two sides. Moreover, in all of these numerical examples, when the probability of failure of the final designs are estimated using a Monte Carlo approach they are found to be very close to the targeted values. In the future, complex failure criteria, such as incorporating fatigue conditions, multi-physics, and so on, will be explored with this method. 
% Moreover, it is observed that in all of our examples the final probability of failure using Monte Carlo is higher than that value calculating subset simulation which is understandable as the number of samples in Monte Carlo are higher than subset simulation. However, the subset simulation and large deviation theory combination can effectively estimate the probability of failure and its gradients without a large number of samples. Although reliability consideration adds more and thicker bars in the designs, the outputs are reliable and optimized.

% Therefore, our method combines subset simulation with large deviation theory and evaluates the probability of failure only every $20$ iterations. In addition, using mini batch stochastic gradient descent allows the optimization effectively happening incorporating uncertainty with a small number of samples. These features make our approach practical for large-scale RBTO problems.
% \appendix 
\section*{Appendix A: Subset Simulation \cite{AU2001loc}} \label{sec:appendix}
Subset simulation was developed to efficiently estimate small failure probabilities which is common in reliability-based analysis of engineering problems \cite{AU2001loc,li2010design}. The main idea behind this method is to decompose a rare failure event into a sequence of more frequent intermediate events. The final probability of failure is the product of probability of failure at each intermediate levels. The conditional probability at these intermediate events are large enough to be estimated using only a small number of samples. Since MC methods are inefficient for estimating small probabilities, subset simulation employs a Markov chain Monte Carlo technique based on the Metropolis-Hastings algorithm to generate samples at each intermediate levels \cite{AU2001loc,li2016matlab,de2021reliability}. These steps are illustrated in Algorithm~\ref{alg:subset}. 

\begin{algorithm}[!htb] 
\caption{Subset simulation for probability of failure assessment \cite{AU2001loc, li2016matlab, de2021reliability}} \label{alg:subset}
\begin{algorithmic}[1] 
\State \textbf{Input:} Density distribution $\rhoo$, number of samples $N_s$ per level, conditional probability $p_0$ (typically $0.1 \le p_0 \le 0.3$) 
\State Generate $N_s$ \textit{i.i.d.} samples $\{\xii_i\}_{i=1}^{N_s}$ of the uncertain variables from pdf $\pdf(\xii)$ 
    \For{$i = 1$ to $N_s$}
        \State Solve $K(\rhoo_k;\xii_i)U(\rhoo_k;\xii_i) = F(\xii_i)$ to obtain displacement $U(\rhoo_k, \xii_i)$
        \State Evaluate the performance function $g(\rhoo_k;\xii_i)$ 
    \EndFor

    \State Sort performance such that
    \[
        g(\boldsymbol{\rho_k}; \boldsymbol{\xi}_1) > g(\boldsymbol{\rho_k}; \boldsymbol{\xi}_2) > \cdots > g(\boldsymbol{\rho_k}; \boldsymbol{\xi}_{N_s})
    \]
        \State Set $b_0 = g(\boldsymbol{\rho_k}; \boldsymbol{\xi}_{\lfloor N_s p_0 \rfloor})$
    \State Set $j = 0$
    
    \While{$b_j > 0$}
    \State Define the intermediate failure criteria
    \[
        \mathcal{F}_j := \{ \boldsymbol{\xi} : g(\boldsymbol{\rho_k}; \boldsymbol{\xi}) \le b_j \}
        \]
    \For{$i = 1, \ldots, \lfloor N_s p_0 \rfloor$}
        \State Generate $\lceil 1/p_0 \rceil$ samples from a Markov chain with stationary pdf $\pdf(\boldsymbol{\xi} \mid \mathcal{F}_j)$
        \State Initialize the chain from the $i$-th sample in $\mathcal{F}_j$
    \EndFor
    \State Sort the new $N_s$ samples as
    \[
        g(\rhoo_k; \boldsymbol{\xi}_1) > g(\rhoo_k; \boldsymbol{\xi}_2) > \cdots > g(\rhoo_k; \boldsymbol{\xi}_{N_s})
    \]
    \State Set $j = j + 1$
    \State Set $b_j = g(\rhoo_k; \boldsymbol{\xi}_{\lfloor N_s p_0 \rfloor})$
    \EndWhile
    
    \State Find the $N_f$ to satisfy
    \[
        g(\rhoo_k; \boldsymbol{\xi}_{N_f}) \ge 0 > g(\rhoo_k; \boldsymbol{\xi}_{N_f+1})
        \]
    \State Estimate the failure probability
    \[
        P_f(\rhoo_k) = \frac{N_f}{N_s} p_0^{\,j}
    \]
\State \textbf{Output:} Probability of failure $P_f$
\end{algorithmic}
\end{algorithm}

\newpage
\bibliographystyle{unsrt}
\bibliography{references}  %%% Uncomment this line and comment out the ``thebibliography'' section below to use the external .bib file (using bibtex) .





% \end{thebibliography}

\end{document}